\documentclass{ieeeaccess}
\usepackage{cite} %new
\usepackage{textcomp} %new
\usepackage{amsmath,amssymb,amsfonts,epsfig,color,graphicx,eucal,mathrsfs}
\usepackage{mathptmx}
\usepackage[cal=boondoxo]{mathalfa}
\usepackage[section]{algorithm}% http://ctan.org/pkg/algorithms
\usepackage{algpseudocode}% http://ctan.org/pkg/algorithmicx
\usepackage{etoolbox}
\usepackage{footnote}
\usepackage{threeparttable}
\usepackage{multirow,multicol}
\usepackage{epstopdf} 
\usepackage{graphicx}
\usepackage{subfigure}  % support sub-figure
\usepackage{caption,setspace}
\usepackage{bigstrut}
\makesavenoteenv{tabular}
\def\BibTeX{{\rm B\kern-.05em{\sc i\kern-.025em b}\kern-.08em
    T\kern-.1667em\lower.7ex\hbox{E}\kern-.125emX}}

\begin{document}
	
\history{}
\doi{}

\title{Reusing Preconditioners in Projection based Model Order Reduction Algorithms}
\author{\uppercase{Navneet Pratap Singh}\authorrefmark{1} AND
%	\IEEEmembership{Fellow, IEEE},
\uppercase{Kapil Ahuja}\authorrefmark{2}
%\IEEEmembership{Member, IEEE}
}
\address{Data \& Computational Sciences Laboratory, Indian Institute of Technology Indore, Indore 453552, India}
%	Computational Science and Engineering Laboratory, IIT Indore, Indore 453552, India}
\address[1]{(e-mail: navneet.diat@gmail.com, phd1301201002@iiti.ac.in)}
\address[2]{(e-mail: kapsahuja22@gmail.com, kahuja@iiti.ac.in)}
\newcommand{\orcidauthorA}{0000-0001-6403-2229} % Add \orcidA{} behind the author's name
\newcommand{\orcidauthorB}{0000-0001-9640-4437} % Add \orcidB{} behind the author's name

\tfootnote{``Further, the second author would like to duly acknowledge the support from MATRICS Scheme of DST-SERB, India (MTR/2017/001023).''}

\markboth
{N. P. Singh, K. Ahuja: Reusing Preconditioners in Projection based Model Order Reduction Algorithms.}
{}%K. Ahuja \headeretal: Preparation of Papers for IEEE TRANSACTIONS and JOURNALS}
\corresp{Corresponding author: Navneet Pratap Singh (e-mail: navneet.diat@gmail.com).}

\begin{abstract}
Dynamical systems are pervasive in almost all engineering and scientific applications. Simulating such systems is computationally very intensive. Hence, Model Order Reduction (MOR) is used to reduce them to a lower dimension. Most of the MOR algorithms require solving large sparse sequences of linear systems. Since using direct methods for solving such systems does not scale well in time with respect to the increase in the input dimension, efficient preconditioned iterative methods are commonly used.
In one of our previous works, we have shown substantial improvements by reusing preconditioners for the parametric MOR (Singh et al. 2019). Here, we had proposed techniques for both,  the non-parametric and the parametric cases, but had  applied them only to the latter. We have four main contributions here. First, we demonstrate that preconditioners can be reused more effectively in the non-parametric case as compared to the parametric one because of the lack of parameters in the former. Second, we show that reusing preconditioners is an art and it needs to be fine-tuned for the underlying MOR algorithm. Third, we describe the pitfalls in the algorithmic implementation of reusing preconditioners. Fourth, and final, we demonstrate this theory on a real life industrial problem (of size 1.2 million), where savings of upto $64 \%$ in the total computation time is obtained by reusing preconditioners. In absolute terms, this leads to a saving of $5$ days. 
\end{abstract}


\titlepgskip=-15pt

\maketitle
%%%%%%%%%%%%%%%%%%%%%%%%%%%%%%%%%%%%%%%%%%%%%%%%%%%%%%%%%%%%%%%%%%%%%%
\section{Introduction}
\label{sec:intro}
{\color{black} Dynamical systems arise in many engineering and scientific applications such as weather prediction, machine design, circuit simulation, biomedical engineering, etc. 
	%There are mainly two classes of dynamical systems; 1) parametric dynamical systems, 2) non-parametric dynamical systems. 
	%In parametric dynamical systems, parameters are linearly embedded in the system matrices. However, non-parametric dynamical systems are free from the parameters.
	Generally, dynamical systems corresponding to real-world applications are extremely large in size. A set of equations describing a parametric nonlinear {second-order} dynamical system is represented as
	\begin{align}\label{eq:second-order-para-nonlinear}
		\begin{split}
			%			\frac{d^2}{dt^2}(g(x(t), \mathfrak{p})) & = \frac{d}{dt}(f(x(t), \mathfrak{p})) + h(x(t), \mathfrak{p}, u(t)), \\ y(t) & = C^T x(t), \\
			g(\ddot{x}(t), \mathfrak{p}) & = f(\dot{x}(t), \mathfrak{p}) + h(x(t), \mathfrak{p}, u(t)), \\ y(t) & = C^T x(t), \\
			%			 \text{given that,}\\
			%			\ddot{x}(t)  & = g(x(t), \mathfrak{p}), \\
			%			\dot{x}(t) & = f(x(t), \mathfrak{p}),
		\end{split}
	\end{align}
	where $t$ is the time variable, $x(t): \mathbb{R} \rightarrow \mathbb{R}^{n}$ is the state, $\mathfrak{p} = (p_1, p_2, \ldots, p_k)$ is the set of parameters (with $p_{\mathcal{j}} \in \mathbb{R}$; for $\mathcal{j} = 1, \ \ldots,  \ k$), $u(t):\mathbb{R} \rightarrow \mathbb{R}^{m}$ is the input, $y(t):\mathbb{R} \rightarrow  \mathbb{R}^{q}$ is the output, $C^T \in \mathbb{R}^{q \times n}$ is the output matrix, and $g(\cdot): \mathbb{R}^{n+k} \rightarrow \mathbb{R}^{n}$, $f(\cdot): \mathbb{R}^{n+k} \rightarrow \mathbb{R}^{n}$ and $h(\cdot): \mathbb{R}^{n+k+m} \rightarrow \mathbb{R}^{n}$ are some nonlinear functions \cite{Ogata2001,  Rewienski2003,  Antoulas2004,  Antoulas05,  Socolar2006,  Bradley2006}. If $m$ and $q$ both are equal to one, then we have a Single-Input Single-Output (SISO) system. Otherwise, it is called a Multi-Input Multi-Output (MIMO) ($m$ and $q > 1)$ system.
	%In (\ref{eq:second-order-para-nonlinear}), if we can linearize these vector-valued functions $\left(\text{i.e.}, \ g(\cdot), \ f(\cdot), \ h(\cdot)\right)$ with respect to the state variable 
	The functions $g(\cdot), \ f(\cdot), \ \text{and} \ h(\cdot)$ are usually simplified as \cite{Rewienski2003, Bradley2006} %, then we have %, then we have %the following linearized state equation and output equation as:   
	\begin{align}\label{eq:linearization_wrt-para}
		\begin{split}
			%D(t;p)\dot{x}(t;p)  &= -K(t;p)x(t;p) + F(t;p)u(t),\\ y(t;p)  &= C(t;p) x(t;p) + G(t;p) u(t),
			g(\ddot{x}(t), \mathfrak{p}) & = \sum_{{\mathcal{j}} = 1}^{k} \mathsf{g}_{\mathcal{j}}(\mathfrak{p}) \mathcal{g}(\ddot{x}(t)), \\
			f(\dot{x}(t), \mathfrak{p}) & = \sum_{{\mathcal{j}} = 1}^{k} \mathsf{f}_{\mathcal{j}}(\mathfrak{p}) \mathcal{f}(\dot{x}(t)),\\
			h(x(t), \mathfrak{p}, u(t)) & = \sum_{{\mathcal{j}} = 1}^{k} \mathsf{h}_{\mathcal{j}}(\mathfrak{p}) \mathcal{h}(x(t),u(t)),
		\end{split}
	\end{align}  
	where $\mathsf{g}_{\mathcal{j}}(\cdot), \ \mathsf{f}_{\mathcal{j}}(\cdot), \ \mathsf{h}_{\mathcal{j}}(\cdot) : \mathbb{R}^{k} \rightarrow \mathbb{R}$  are scalar-valued functions while $\mathcal{g}(\cdot), \ \mathcal{f}(\cdot) : \mathbb{R}^{n} \rightarrow \mathbb{R}^{n}$, and $\mathcal{h}(\cdot): \mathbb{R}^{n+m} \rightarrow \mathbb{R}^{n}$ are vector-valued. %can be nonlinear functions. %In (\ref{eq:linearization_wrt-para}), 
	Next, we look at simplifications to (\ref{eq:second-order-para-nonlinear}) based upon the three predicates; the presence of parameters; the degree of non-linearity, and the order of the system.    
	\begin{itemize}
		\item If $\mathsf{g}_{\mathcal{j}}(\mathfrak{p}), \ \mathsf{f}_{\mathcal{j}}(\mathfrak{p}), \ \text{and} \ \mathsf{h}_{\mathcal{j}}(\mathfrak{p})$ are independent of the parameters, then (\ref{eq:second-order-para-nonlinear}) becomes a non-parametric dynamical system.
		\item Bilinear systems are one of the common types of nonlinear dynamical systems. Here, there is a product between the state variables and the input variables.	
		Another important class of nonlinear dynamical systems is the quadratic systems. Here, there is product among the state variables. 	If $\mathcal{g}(\cdot)$ and $\mathcal{f}(\cdot)$ are linear functions of the state variables, and $\mathcal{h}(\cdot)$ is a linear function of the state and the input variables, then (\ref{eq:second-order-para-nonlinear}) is called a linear dynamical system.	
		\item Finally, if the second derivative term in (\ref{eq:second-order-para-nonlinear}) is not present, then (\ref{eq:second-order-para-nonlinear}) becomes a {first-order} dynamical system.	
	\end{itemize}	 
	%	In nonlinear dynamical system, bilinear dynamical system are commonly used, here, there is a product between state and input variables. There is one more class in the bilinear dynamical system, that is called quadratic-bilinear dynamical system, here, there is square-term with state variable. This class is also commonly studied.
	%The in nonlinear system, bilinear is commonly used, there is product between state variable $x(t)$ and input variable $u(t)$. And there is one more case of bilinear, that quadratic bilinear, where there is square term with state variable.  There variation of these for example quadratic bilinear system is also commonly studied.
	%If $\mathcal{g}(\cdot), \ \mathcal{f}(\cdot)$ are linear function of $x$, and $\mathcal{h}(\cdot)$ is linear function of $x$ and $u$, then (\ref{eq:second-order-para-nonlinear}) becomes a linear dynamical system. 
	%Otherwise it is represented as nonlinear dynamical system. 
	
}

{\color{black}Simulation of large dynamical systems can be unmanageable due to high demands on computational resources. These large systems can be reduced into a smaller dimension by using Model Order Reduction (MOR) techniques \cite{Grimme1997,  Gugercin2003,  Antoulas05,  Schilders2008, Gugercin2008, Breiten2013}. The reduced system has approximately  the same characteristics as the original system but it requires significantly less computational effort in simulation. MOR can be done in many ways such as balanced truncation, Hankel approximations, and Krylov projection \cite{Grimme1997, Gugercin2003, Antoulas05, Breiten2013}. Among these, the projection methods are quite popular, and hence, we focus on them.%, and the reduced system can be obtained by projecting the original system on a lower dimensional subspace.
	
Some of the commonly used projection-based MOR algorithms for different types of dynamical systems are summarized in Table \ref{tab:MOR-algo}.
	% \usepackage{multirow}
	      \begin{table*}[]
      	%\normalsize
		\centering
		\caption{MOR Algorithms.}
		\setlength{\tabcolsep}{3pt}
	%	\medskip % not necessary if you load caption
%		\def\arraystretch{1.0}
		\begin{tabular}{|c|l|l|l|c|}%{p{0.25\linewidth}p{0.25\linewidth}p{0.25\linewidth}}
			%Content\footnotemark\\
			\hline
			{\multirow{2}{*}{S. No.}} & \multicolumn{2}{c|}{\multirow{2}{*}{{Category}}} & \multicolumn{2}{c|}{{Order}} \\ \cline{4-5} 
			& \multicolumn{2}{c|}{} & \multicolumn{1}{c|}{{First}} & {Second} \\ \hline
			1. & \multicolumn{2}{l|}{Non-parametric Linear} & IRKA \cite{Gugercin2008}, $(Sy)^{2}\text{IRKA}$ \cite{Breiten2013} & \multicolumn{1}{l|}{\begin{tabular}[c]{@{}l@{}}SOR-IRKA \cite{Wyatt2012}, SO-IRKA \cite{Qiu2018} \\ SOSPDR \cite{BaiS05}, AIRGA \cite{Bonin20161} \end{tabular}} \\ \hline
			{\multirow{2}{*}{2.}} & \multicolumn{1}{c|}{\multirow{2}{*}{Non-parametrc}} & Bilinear & BIRKA \cite{Benner2012}, TB-IRKA \cite{Flagg2012, Choudhary2019} & -- \\ \cline{3-5} 
			& \multicolumn{1}{c|}{} & Quadratic-bilinear & QB-IHOMM \cite{Mian2019} & -- \\ \hline
			3. & \multicolumn{2}{l|}{ Parametric Linear} & I-PMOR \cite{Baur2011}, RPMOR \cite{Benner2014} & \multicolumn{1}{l|}{RPMOR \cite{Feng2013}}  \\ \hline
			{\multirow{2}{*}{4.}} & \multirow{2}{*}{Parametric} & Bilinear & I-PMOR-Bilinear \cite{Rodriguez2018} & -- \\ \cline{3-5}  
			& & Quadratic-bilinear & QB-IRKA \cite{Cao2019} & -- \\ \hline
		\end{tabular}
		%	\\Content\footnotemark
		%	\footnote{Iterative Rational Krylov Algorithm (IRKA), IRKA for symmetric Sylvester equations $((Sy)^{2}\text{IRKA})$, Second-Order Structure-Preserving Dimension Reduction (SOSPDR), Adaptive Iterative Rational Global Arnoldi (AIRGA),  Bilinear IRKA (BIRKA), Truncated Bilinear IRKA (TB-IRKA), Implicit Higher-Order Moment Matching technique (QB-IHOMM), PMOR with Interpolatory Projections (I-PMOR), Robust Algorithm for Parametric Model Order Reduction (RPMOR), Interpolatory Model Reduction of Parameterized Bilinear dynamical systems (I-PMOR-Bilinear), Quadratic Bilinear Iterative Rational Krylov algorithm (QB-IRKA).}
		%Content\footnotetext{}
		\label{tab:MOR-algo}
	\end{table*}

}

{\color{black}In the above mentioned MOR algorithms,  sequences of very large and sparse linear systems arise during the model reduction process. Solving such linear systems is the main computational bottleneck in efficient scaling of these MOR algorithms for reducing extremely large dynamical systems.  
	Preconditioned iterative methods are commonly used for solving such linear systems \cite{Ahuja2012, Ahuja2015}. In most of the above listed MOR algorithms, the change from one linear system to the next is usually very small, and hence, the applied preconditioner could be reused.  

	Next, we briefly summarize the past work that has been done in the field of reusing preconditioners. 
	%in MOR context and others context also. Although, our focus is only on a non-parametric dynamical system ($1^{st}$ and $2^{nd}$ categories as listed above). {\it First-time} 
	This technique was first applied in the QMC context, where it was referred to as recycling preconditioner \cite{ahujathesisPhD, ahuja2011improved}. In the optimization context, this approach was applied in \cite{McNally2017}, where it was termed as the preconditioner update. Such a technique was first applied in the MOR context in \cite{Wyatt2012}, and more recently in \cite{grim2015reusing}, where the focus was mostly on MOR of non-parametric linear {first-order} dynamical systems (part of the first category above).

	%	it has been used as a preconditioner update \cite{McNally2017}. One of the first works has been done in the MOR context \cite{Wyatt2012}, for the non-parametric linear {\it first-order} dynamical system with the IRKA algorithm. Further, this work has been extended for other applications of the dynamical system \cite{grim2015reusing}. 
	%------------
	The main goal of this paper is to demonstrate the reuse of preconditioners in the remainder of the algorithms for the first category above (MOR of non-parametric linear {second-order} dynamical systems) as well as the algorithms for the second category above (MOR of non-parametric bilinear/ bilinear-quadratic dynamical systems). %In general, the focus here is, on MOR of dynamical systems that are non-parametric.
	
	In one of our recent works \cite{navneet2019}, we had proposed a general framework for reuse of preconditioners during MOR of both non-parametric and parametric dynamical systems. However, in \cite{navneet2019} we had demonstrated application of this framework for the parametric case only. That is, the third category above (MOR of parametric linear dynamical systems). We are currently (and separately) working on the algorithms for the fourth category above as well (MOR of parametric bilinear/ bilinear-quadratic dynamical systems).

	To summarize, in this paper we broadly demonstrate the application of our above mentioned framework for MOR of non-parametric dynamical systems. We have four contributions as below, which have not been catered in any of the above cited papers. 

	{\color{black} %There are four important innovative aspects of this work, which have not been catered in any of the above-cited papers.
		\begin{itemize}
			\item[(i)] We demonstrate that the reuse of preconditioners can be done more effectively in the non-parametric case as compared to the parametric case because of the lack of parameters in the former.   
			
			\item[(ii)] We show that as the underlying MOR algorithms get more intricate, the reuse of preconditioners needs to be fine-tuned. 
			%As mentioned above, we focus on following three algorithms: AIRGA algorithm (discuss in the next section), optimal $H_2$ model reduction method (discuss in Appendix I), and Nonlinear Parameterized Model Order Reduction (NLPMOR) algorithm (discuss in Appendix II) where the reuse preconditioner closely follows the $2^{nd}$ approach of the frameworks \cite{navneet2019}, as discussed above.
			\item[(iii)] We highlight that there are multiple pitfalls in the algorithmic implementation of 
			reusing preconditioners, which if not done efficiently, could actually increase the computational complexity of the underlying MOR algorithms instead of reducing it.
			%	doesn't increase the computational complexity by applying this technique to underlying algorithms.  
			
			\item[(iv)] We experiment on a massively large and real-life industrial problem (BMW disc brake model), which is of size $1.2$ million. %unknowns, by using the AIRGA algorithm. 
			We reduce the total computation time from $197$ hours to about $72$ hours (approximately), leading to a saving of $64 \%$.

		\end{itemize}
		
		The paper has four more sections. We discuss MOR techniques in Section \ref{ses:MOR}. The theory of reusing preconditioners is described in Section \ref{sec:Reusing_Precond.}. We support our theory with numerical experiments in Section \ref{sec:Numerical_results}. Finally, conclusions and future works are discussed in Section \ref{sec:conl_fut}. For the rest of this paper, $\|\cdot\|_{f}$ denotes the Frobenius norm,  $\|\cdot\|$  denotes the Euclidean norm for vectors and the induced spectral norm for matrices, $\otimes$ refers to the Kronecker product (i.e. an operation on two matrices of arbitrary size), $vec(\cdot)$ signifies the vectorization of a matrix, and $I$ denotes the Identity matrix.}

%%%%%%%%%%%%%%%%%%%%%%%%%%%%%%%%%%%%%%%%%%%%%%%%%%%%%%%%%%%%
\section{MOR}
\label{ses:MOR}
{\color{black}

	{\color{black}
		%We have briefly discussed the details of Ritz-Galerkin projection based with one of the MOR algorithm, i.e., AIRGA in Algorithm \ref{ALGO:AIRGA}. 	
		%Again from the listed projection-based MOR algorithms in the previous section, out of them, IPMOR, IRKA, $(Sy)^2$IRKA, IRKA for MIMO fall in this class.
	}
	As above, our focus is on MOR of the non-parametric dynamical systems. Hence, we summarize some of the previously listed such algorithms here. 	
	AIRGA \cite{Bonin20161} is a Ritz-Galerkin projection based algorithm for MOR of linear {second-order} MIMO dynamical systems with proportional damping, which for the MIMO case are represented as %SOSPDR is mainly designed for Single-Input Single-Output (SISO) systems. Here, the standard Arnoldi method for the a second-order Krylov subspace (which is referred to
	%as second-order Arnoldi method or SOAR method) is used for generating the reduced systems matrices.
	%AIRGA is mainly design for Multi-Input Multi-Output (MIMO) dynamical systems with proportional damping.
	%is mainly proposed for second-order non-parametric Multi-Input Multi-Output (MIMO) dynamical systems with proportional damping. 
	%	A set of equations describing a non-parametric {\it second-order} MIMO dynamical system is represented as
	\begin{align}\label{eq:linear-non-para}
	\begin{split}
	M \ddot{x}(t) & = - D \dot{x}(t) - K x(t) + F u(t), \\
	y(t) & =  C^T x(t),
	\end{split}  
	\end{align}        	   	 
	where $M,\ D,\ K \in  \mathbb{R}^{n \times n}, \ F \in  \mathbb{R}^{n \times m}, \ C \in  \mathbb{R}^{n \times q},$ and $D = \alpha M + \beta K$. Here, $\alpha, \ \beta$ are some scalar values. 
	%The projection method is used for constructing the above equations into a lower dimension subspace, 
	Let $V \in \mathbb{R}^{n \times r}$ and its columns span a $r$-dimension subspace ($r \ll n$).	
	In principle, the Ritz-Galerkin projection method involves the steps below.
	
	%Let the state variable $x(t)$ is projected onto a smaller dimensional subspace. Let $V \in \mathbb{R}^{n \times r}$ is orthogonal matrix with orthonormal columns and using $x(t) \approx V \hat{x}(t)$ in the system (\ref{eq:linear-non-para}), we obtain the following:  
	\begin{itemize}
		\item Approximating the reduced state vector $\hat{x}(t)$ using ${V}$ as $ x(t) \approx V \hat{x}(t)$ leads to
		\begin{align*}%\label{TD_r11}
		\begin{split}
		M V\ddot{\hat{x}}(t) + D V \dot{\hat{x}}(t)+KV \hat{x}(t) - F u(t) & = r(t),\\
		\hat{y}(t) & = \ C^T V \hat{x}(t),
		\end{split}
		\end{align*}
		where $r(t)$ is the residual after projection. %Applying the Ritz-Galerkin approach by multiplying $V^T$ in the first equation of (\ref{TD_r11}), we get 
		\item Enforcing the residual $r(t)$ to be orthogonal to $V$ or $V^T r(t) = 0$ %. That is,
		%		\begin{align*}%\label{TD_r1}
		%		\begin{split}
		%		V^{T}\left(M V\ddot{\hat{x}}(t) + D V \dot{\hat{x}}(t)+KV \hat{x}(t) - F 
		%		u(t)\right) & = 0,\\
		%		\hat{y}(t) & =  C^T V \hat{x}(t). %+C_{v}V\dot{\hat{x}}(t),
		%		\end{split}
		%		\end{align*}
		%		\item This 
		leads to the reduced system  given as follows:
		%Comparing (\ref{TD_r1}) with (\ref{TD_r}) yields
		\begin{align*}%\label{eq:red_sys}
		\begin{split}
		\hat{M} \ddot{\hat{x}}(t) + \hat{D} \dot{\hat{x}}(t)+ \hat{K} \hat{x}(t) - \hat{F}	u(t) & = 0,\\
		\hat{y}(t) & =  \hat{C}^T  \hat{x}(t),
		\end{split}
		\end{align*}	
	\end{itemize}
	where $\hat{M} = V^{T}MV,\ \hat{D} = V^{T}DV,\ \hat{K}=V^{T}KV,\ \hat{F} = V^{T}F, \ \text{and} \ \hat{C}^T = C^TV$. To compute this projection matrix $V$, AIRGA matches the moments of the original system transfer function and the reduced system transfer function.
	We briefly summarize AIRGA in Algorithm \ref{ALGO:AIRGA}, where parts relevant to solving linear systems are only listed.	
	\begin{algorithm}[]
		%\large
		%\normalsize
		%\setstretch{1.2}
		\caption{:~AIRGA~\cite{Bonin20161}}
		\label{ALGO:AIRGA}
		\begin{algorithmic}[1]
			%\While {$||\hat{H}_{old}-\hat{H}_{new}||_{H_{2}} \ \leq \ tol$}
			\Statex  Input: {$M,\ D,\ K,\ F,\ C$; \ $S$ is the set of initial expansion points $s_{i},\ i=1,\ \ldots, \ \ell$}.
			\Statex Output: ${\hat{M}, \ \hat{D}, \ \hat{K}, \ \hat{F}, \ \hat{C}}$.
			\State $z = 1$
			\While{(no convergence)}
			\For {$i = 1,\ \ldots,\ \ell$}  
			
			\State  $X^{(0)}(s_{i})= (s_{i}^{2}M+s_{i}D+K)^{-1}F$	
			\Statex
			\State 	$V_1 = \frac{X^{(0)}(s_{i})}{\|X^{(0)}(s_{i})\|_f}$ 
			\EndFor
			%\While {$||\hat{H}_{j}-\hat{H}_{j+1}||_{H_{2}} \leq tol $}
			\State j = 1
			\While{(no convergence)} %and $j \le \lceil r_\text{max} / m \rceil $}
			\For {$i = 1,\ \ldots, \ \ell$}
			\State  $X^{(j)}(s_{i})= -(s_{i}^{2}M+s_{i}D+K)^{-1}MV_j$
			\Statex	
			\State 	$V_{j+1} = \frac{X^{(j)}(s_{i})}{\|X^{(j)}(s_{i})\|_f}$ %for $i = \{1,\ldots,l\}$
			\EndFor
			%\State 	$V_{j+1} = \frac{X^{(j)}(s_{i})}{\|X^{(j)}(s_{i})\|_f}$ for $i = \{1,\ldots,l\}$
			\State $j = j+1$
			\EndWhile
			\State  \text{{\it ``All the given set of expansion points}} 					
			\Statex \qquad \text{{(i.e. $s_1, \ s_2,  \ \ldots, \ s_{\ell}$) are updated''}}
			\State $z = z+1$
			\EndWhile
			\State $\hat{M} = V^{T}MV,\ \hat{D} = V^{T}DV,\ \hat{K}=V^{T}KV,\ \hat{F} = V^{T}F, \ \text{and} \ \hat{C}^T = C^TV$
		\end{algorithmic}
	\end{algorithm}
	
	BIRKA \cite{Benner2012} is a Petrov-Galerkin projection based algorithm for MOR of the {bilinear} {first-order} dynamical systems, which for the MIMO case are represented as 
	%A set of equations describing a non-parametric {bilinear}  {\it first-order} MIMO dynamical system is represented as
	%A non-parametric bilinear {\it first-order} MIMO dynamical system can be represented as
	\begin{align}\label{eq:non-para-bilinear}
	\begin{split}
	\dot{x}(t) & = Kx(t) + \sum_{\mathsf{j}=1}^{m} N_{\mathsf{j}}x(t)\mathsf{u}_{\mathsf{j}}(t)  + F u(t), \\
	y(t) & = C^Tx(t),
	\end{split}
	\end{align}
	where $K, \ N_{\mathsf{j}} \in \mathbb{R}^{n \times n}, \ F \in \mathbb{R}^{n \times m}, \ C \in \mathbb{R}^{n \times q}$, and $u = [\mathsf{u}_1, \ \mathsf{u}_2, \ \ldots, \mathsf{u}_m] \in \mathbb{R}^m$.
	%We project the above equation into a lower dimension subspace, let $V \in \mathbb{R}^{n \times r}, \text{and} \ W \in \mathbb{R}^{n \times r}$ is orthogonal matrix with orthonormal columns and using $x(t) \approx V \hat{x}(t)$ in the system (\ref{eq:non-para-bilinear}), we obtain the following by enforce the Petrov-Galerkin condition
	%The projection method is used for constructing the above equations into a lower dimension subspace, 
	Let columns of  $V, W \in \mathbb{R}^{n \times r}$ span two $r$-dimension subspaces (where, as earlier, $r \ll n$ ).  In principle, the Petrov-Galerkin projection method involves the steps below.
	\begin{itemize}
		\item Approximating the reduced state vector $\hat{x}(t)$ using  ${V}$ as $ x(t) \approx V \hat{x}(t)$ leads to
		\begin{align*}%\label{eq:two_proj1}
		\begin{split}
		V\dot{\hat{x}}(t) - KV\hat{x}(t) - \sum_{\mathsf{j}=1}^{m} N_{\mathsf{j}}V\hat{x}(t)\mathsf{u}_{\mathsf{j}}(t)  - F u(t) & = r(t),\\
		\hat{y}(t) & = C^T V \hat{x}(t),
		\end{split}
		\end{align*}
		where $r(t)$ is the residual after projection.
		%i.e., %the state variable $x(t)$ is projected onto a smaller dimensional subspace. %Let $W \in \mathbb{R}^{n \times r}$ and $V \in \mathbb{R}^{n \times r}$ are two orthogonal matrices with orthonormal columns and using 
		%$x(t) \approx V \hat{x}(t)$ and enforce the Petrov-Galerkin condition
		%\begin{align}\label{eq:two_proj1}
		%\begin{split}
		%& M V\ddot{\hat{x}}(t) + D V \dot{\hat{x}}(t)+KV \hat{x}(t) - F 
		%u(t)= r(t),\\
		%& \hat{y}(t)= \ C_{p} V \hat{x}(t)+C_{v}V\dot{\hat{x}}(t).
		%\end{split}
		%\end{align}
		%%where $r(t)$ is residual after projection. 
		%Applying the Petrov-Galekin approach by multiplying $W^T$ in the first equation of (\ref{eq:two_proj1}) we get
		%let $W \in \mathbb{R}^{n \times r}$ and $V \in \mathbb{R}^{n \times r}$ are the two orthogonal matrices and the system (\ref{eq1}) can be projected as  	
		\item Enforcing the residual $r(t)$ to be orthogonal to $W$ or $W^T r(t) = 0$ %. That is,
		%		\begin{align*}%\label{eq:two_proj}
		%		\begin{split}
		%		W^{T}\left(V\dot{\hat{x}}(t) - KV\hat{x}(t) - \sum_{\mathsf{j}=1}^{m} N_{\mathsf{j}}V\hat{x}(t)\mathsf{u}_{\mathsf{j}}(t)  - F u(t)\right) & = 0,\\
		%		\hat{y}(t) & =  C^T V \hat{x}(t),
		%		\end{split}
		%		\end{align*}
		%	
		%	\item This 
		leads to the reduced system given by
		%On comparing (\ref{eq:two_proj}) with (\ref{TD_r})
		\begin{align*}%\label{eq:red_sys1}
		\begin{split}
		\dot{\hat{x}}(t) - \hat{K}\hat{x}(t) - \sum_{\mathsf{j}=1}^{m} \hat{N}_{\mathsf{j}}\hat{x}(t)\mathsf{u}_{\mathsf{j}}(t)  - \hat{F} u(t) & = 0,\\
		\hat{y}(t) & = \ \hat{C}^T  \hat{x}(t),
		\end{split}
		\end{align*}
	\end{itemize}
	where $\hat{K} = (W^{T}V)^{-1}W^{T}KV,\ \hat{N}_{\mathsf{j}} = (W^{T}V)^{-1} W^{T}N_{\mathsf{j}}V, \ \hat{F} = (W^{T}V)^{-1} W^{T}F, \ \hat{C}^T= C^TV$, and $(W^{T}V)^{-1}$ is assumed to be invertible. Here, $V$ and $W$ are computed by using interpolation, where the original system transfer function and its derivative are respectively matched with the reduced system transfer function and its derivative at a set of points. 
	%The algorithm is iterative in nature and for a good initial guess of the reduced system, converges to the ideal interpolation points. 
	We briefly summarize BIRKA in Algorithm \ref{Algo:BIRKA}, where again, only parts related to solving linear systems are listed. 
	%	We have listed this in Algorithm \ref{Algo:BIRKA}.
	%Different solutions of subspace the $V$ and $W$ give different reduced systems, but the choices of the subspace $V$ and $W$ by enforcing interpolation. There are two different ways of doing the interpolation, i.e., subsystem interpolation and Volterra series interpolation \cite{Benner2012}. BIRKA is based on Volterra series interpolation, and it gives a locally $H_2$-optimal reduced system. 
	\begin{algorithm}[]
			%	\setstretch{1.2}
		\caption{:~BIRKA~\cite{Benner2012}}
		\label{Algo:BIRKA}
		{\color{black}\begin{algorithmic}[1]
				\Statex Input $K,\ N_1, \ \ldots, \ N_m , \ F,\ C$, and initial guess of the reduced system $\check{K},\ \check{N}_1, \ \ldots, \ \check{N}_m,\ \check{F},\ \check{C}$ %Also select stopping tolerance $btol$.
				\Statex Output $\hat{K}, \ \hat{N}_{1}, \ \ldots, \ \hat{N}_{m} , \ \hat{F}, \ \text{and} \ \hat{C}$
				\State $z=1$
				\While{(no convergence)} 
				%\begin{enumerate}
				\State	 $R \Lambda R^{-1} = \check{K},\ \check{\check{F}}= \check{F}^T R^{-T}, \ \check{\check{C}}=\check{C}R, \ \check{\check{N}}_{\mathsf{j}}= R^T \check{N}_{\mathsf{j}} R^{-T}$ \Statex \quad  \ \ for $\mathsf{j}=1 ,\ \ldots,\ m$ 
				\Statex
				\State	 $vec\left ( V \right ) = \left ( -\Lambda \otimes I_n - I_{r} \otimes K - \sum\limits_{\mathsf{j}=1}^{m} \check{\check{N}}^T_{\mathsf{j}} \otimes  N_{\mathsf{j}}   \right )^{-1}$  
				\Statex \qquad \qquad \qquad $\left ( \check{\check{F}}^T \otimes F\right ) \ {vec(I_{m})} $
				\State	 $vec\left ( W \right ) = \left ( -\Lambda \otimes I_n - I_{r} \otimes K^T - \sum\limits_{{\mathsf{j}}=1}^{m}\check{\check{N}}_{\mathsf{j}} \otimes N^T_{\mathsf{j}} \right )^{-1}$ \Statex \qquad \qquad \qquad $\left ( \check{\check{C}}^T \otimes C^T\right ) \ {vec(I_{q})}$
				\Statex
				\State	$V = orth\left ( V  \right ) , \ W = orth\left ( W \right ) $
				\Statex
				\State	 $\check{K}= (W^T V)^{-1} W^T K V$, \   $\check{N}_{\mathsf{j}}= \left ( W^T V \right )^{-1} W^T N_{\mathsf{j}} V,$ \
				\Statex \quad \ $ \check{F}=\left ( W^T V \right )^{-1} W^T F,$ \ 
				$ \check{C}= CV$	
				%\end{enumerate}
				\State $z = z+1$
				\EndWhile
				\State $\hat{K} = \check{K}, \  \hat{N}_{\mathsf{j}} = \check{N}_{\mathsf{j}}, \  \hat{F} =\check{F}, \ \text{and} \ \hat{C} = \check{C}$
				%\Statex 
			\end{algorithmic}}
		\end{algorithm}
		
		QB-IHOMM algorithm \cite{Mian2019} is a Petrov-Galerkin projection based algorithm for MOR of the {quadratic-bilinear} dynamical systems, which for the SISO case are represented as \footnote{ A variant of BIRKA for MOR of the quadratic-bilinear dynamical systems also exists. Preconditioned iterative solves and reusing preconditioners can be applied here as done for BIRKA.
			%	where would be similar as done for BIRKA. 
			Hence, we focus on the QB-IHOMM algorithm that has been developed for the SISO case only.}
		%A set of equations describing a non-parametric {quadratic-bilinear} {\it first-order} SISO dynamical system is represented as
		\begin{align}\label{eq:non-para-Quard-bilinear}
		\begin{split}
		D\dot{x}(t) & = Kx(t) +  Nx(t)u(t) + H\left(x(t) \otimes x(t)\right) + F u(t), \\
		y(t) & = C^Tx(t),
		\end{split}
		\end{align}
		where $D, \ K, \ N \in \mathbb{R}^{n \times n}, \ H \in \mathbb{R}^{n \times n^2}, \ F \in \mathbb{R}^{n \times 1}, \ C \in \mathbb{R}^{n \times 1}$. 
		%Also, $x(t) \in \mathbb{R}^n$ is state, $u(t), \ y(t) \in \mathbb{R}$ are input and output, respectively. 
		Let columns of  $V, \ W \in \mathbb{R}^{n \times r}$ span two $r$-dimension subspaces (where as earlier, $r \ll n$ ).  In principle, the Petrov-Galerkin projection method involves the steps below.
		
		%	Projection method is used for constructing the above equations into a lower dimension subspace, where two matrices $V$ and $W$ are identified such that their columns span $r$-dimension subspace $\mathcal{V}$ and $\mathcal{W}$, respectively.  In principle, projection method involves the following steps:
		\begin{itemize}
			\item As before, approximating the reduced state vector $\hat{x}(t)$ using  ${V}$ as $ x(t) \approx V \hat{x}(t)$ leads to
			\begin{align*}
%			\resizebox{.45\textwidth}{!}
%			{$\begin{aligned}
			%\text{such that} \ \  x(t) \approx V \hat{x}(t), \\
			D V \dot{\hat{x}}(t) - KV\hat{x}(t) & -  NV\hat{x}(t)u(t) \\ & - H\left(V\hat{x}(t) \otimes V\hat{x}(t) \right) - F u(t)  = r(t), \\
			& \qquad \qquad \qquad \qquad  y(t)   = C^T V \hat{x}(t),
%			\end{aligned}$
%		    }
			\end{align*}
			where $r(t)$ is the residual after projection.
			\item Enforcing the residual $r(t)$ to be orthogonal to $W$ or $W^T r(t) = 0$ %. That is,
			%		\begin{align*}
			%		W^T \left( D V \dot{\hat{x}}(t) - KV\hat{x}(t) -  NV\hat{x}(t)u(t) - H\left(V\hat{x}(t) \otimes V\hat{x}(t)\right) - F u(t) \right) & = 0, \\
			%		y(t) & = C^TV\hat{x}(t),
			%		\end{align*}
			%		\item This 
			leads to the reduced system  given by
			\begin{align*}
%			\resizebox{.45\textwidth}{!}
%			{$\begin{aligned}
			\hat{D} \dot{\hat{x}}(t) - \hat{K}\hat{x}(t) -  \hat{N} \hat{x}(t)u(t) - \hat{H}\left(\hat{x}(t) \otimes \hat{x}(t)\right) - \hat{F} u(t)   = 0, \\
			y(t)   = \hat{C}^T\hat{x}(t),
%			\end{aligned}$
%		    }
			\end{align*}
			%\begin{align*}
		\end{itemize}
		where $\hat{D} = W^{T}DV, \ \hat{K} = W^{T}KV, \ \hat{N} = W^{T}NV,$ 
		
		$\hat{H} = W^{T}H(V \otimes V), \ \hat{F} = W^{T}F, \ \hat{C}^T = C^TV.$ Here, $V$ and $W$ are computed by matching the moments of the original system transfer function and the reduced system transfer function.
		We briefly summarize QB-IHOMM in Algorithm \ref{Algo:QB-IHOMM}, where as earlier, only parts related to solving linear systems are listed.
		%We have listed this in Algorithm \ref{Algo:QB-IHOMM}.
		\begin{algorithm}[]
		%	\setstretch{1.2}
			%	\large
			%	\normalsize
			\caption{:~QB-IHOMM~\cite{Mian2019}}
			\label{Algo:QB-IHOMM}
			{\color{black}\begin{algorithmic}[1]
					\Statex Input:  $D,\ K,\ N, \ H, \ F, \ C$; \ interpolation points \ $\sigma_i \in \mathbb{C}$ for $i=1, \ \ldots, \ \ell$; higher orders moments numbers $P, Q \in \mathbb{N}$
					\Statex Output: ${\hat{D}, \ \hat{K}, \  \hat{N}, \ \hat{H}, \ \hat{F}, \ \hat{C}}$
					\State $V = \left[ \ \right], \ \ W = \left[ \ \right]$
					\For {$j = 0, \ \ldots, \ P+Q$} 
					\For {$ i = 1, \ \ldots, \ \ell$}
					\State $X_j(\sigma_i) = [(\sigma_iD-K)^{-1}D]^{j}(\sigma_iD-K)^{-1}F$
					\Statex
					\State $V = \left[V \ \ X_j(\sigma_i) \right]$
					\EndFor    
					\EndFor 
					\For {$j = 0, \ \ldots, \ Q$} 
					\For {$ i = 1, \ \ldots, \ \ell$}
					\State $X_j(2\sigma_i)^T = [(2\sigma_iD-K)^{-T}D^T]^{j}(2\sigma_iD-K)^{-T}C^T$
					\Statex
					\State $W = \left[W \ \ X_j(2\sigma_i)^T \right]$
					\EndFor    
					\EndFor 
					\State $U = orth([V \ W])$    	
					\State Construct the reduced system as
					\Statex $\hat{D} = U^{T}DU, \ \hat{K}=U^{T}KU, \ \hat{N} = U^{T}NU,$ 
					\Statex $\hat{H} = U^{T}H(U \otimes U), \ \hat{F} = U^{T}F, \ \hat{C}^T = C^TU.$			
				\end{algorithmic}}
			\end{algorithm}   
		}
\section{Proposed Work} \label{sec:Reusing_Precond.}
{\color{black}Here, we discuss preconditioned iterative methods in Section \ref{sec:prec-itr}. In Section \ref{sec:thr-resuse}, we revisit the theory of reusing preconditioners from \cite{navneet2019}. Finally, we discuss application of reusing preconditioners to the earlier discussed algorithms in  Section \ref{sec:app-reuse}}.    

%	Our goal here is not only to make the linear system easier to solve by an iterative method. But also, to develop a preconditioner that should be cheap to construct and apply. 
%	%Some existing preconditioning techniques include Successive Over Relaxation, Polynomial, Incomplete Factorizations (ILU), Sparse Approximate Inverse~(SPAI), and Algebraic Multi-Grid~\cite{Benzi2002,benzi1999comparative,Chow1998}. 
%	Incomplete Factorizations (ILU) and Sparse Approximate Inverse~(SPAI) are the most general type of preconditioners. We use 
%	%We use SPAI preconditioner here since these (along with incomplete factorizations) are known to work in the most general setting. 
%	SPAI as a preconditioner since it can be easily parallelized \cite{grote1997}. However, Incomplete Factorizations is difficult to  parallelized. Hence, we use a parallel version of SPAI. Next, we briefly discuss SPAI preconditioner in the contexts of AIRGA algorithm.
\subsection{Preconditioned iterative methods}\label{sec:prec-itr}
%{\color{blue}
%	The MOR algorithms discussed above require solving the sequence of linear systems, which is the main challenging task for reducing the sizeable dynamical system. Here, the linear system matrices are very large and sparse. The direct methods, which are based upon different matrix factorizations, are commonly used for solving these linear systems of equations. However, they become prohibitively expensive for extremely large sizes (hundreds of millions of equations to billions of equations;  as here). In such cases, using iterative methods are usually the only viable option, which scale well both in time and memory. Although iterative methods are not as robust or reliable as direct methods, they are still preferred because scaling is a bigger issue.} 

{\color{black}	Krylov subspace based methods are very popular class of iterative methods \cite{Saad2003, Shen2019}. 
	%There are many types of Krylov subspace methods. Some commonly used ones are Conjugate Gradient (CG), Generalized Conjugate Residual Orthogonal (GCRO), Generalized Minimal Residual (GMRES), Minimum Residual (MINRES), and BiConjugate Gradient (BiCG) etc. 
	%The choice of method is problem dependent. 
	Let 
	$%\begin{align*}
	Ax = b
	$ %\end{align*}
	be a linear system, with $A \in \mathbb{R}^{n \times n}, \ b \in \mathbb{R}^n$, $x_0$ the initial solution and $r_0$ (where $r_0 = b-Ax_0$) the initial residual.
	We find the solution of a linear system in $\mathbb{K}_{\mathcal{k}}(A, \ r_0) = span\{r_0, \ Ar_0, \ A^2r_0, \ \ldots, \ A^{{\mathcal{k}}-1}r_0\}$, where $\mathbb{K}_{\mathcal{k}}(\cdot, \ \cdot)$ represents the Krylov subspace.
	
	%Then, solution of this system can be computed by span of the Krylov subspace. 
	Often iterative methods are slow or fail to converge, and hence, preconditioning is used to accelerate them. If $P$ is a non-singular matrix that approximates the inverse of $A \     (\text{that is}, \ P \approx A^{-1})$, then the preconditioned system becomes $AP\tilde{x} = b$ with $x = P\tilde{x}$ \footnote{This is right preconditioning. Similarly,  center and left preconditioning can be applied \cite{Benzi2002}.}. We expect that the preconditioned iterative solves would find a solution in less amount of time as compared to the unpreconditioned ones. For most of the input dynamical systems (as mentioned here), the Krylov subspace methods fail to converge (see Numerical Experiments section). Hence, we use a preconditioner. 
	
	The goal is to find a preconditioner that is cheap to compute as well as apply. There exist many preconditioning techniques \cite{grote1997, Benzi2002, Chow1998, Fallah2020}, like incomplete factorizations, Sparse Approximate Inverse (SPAI) etc. SPAI preconditioners are known to work in the most general setting and can be easily parallelized. Hence,  we use them.
	
	%, which are known to work in the most general setting. However, SPAI preconditioners can be easily parallelized. Hence, we use a parallel version of SPAI. %We briefly discuss SPAI preconditioner next.
	
	For constructing a preconditioner $P$ corresponding to a coefficient matrix $A$, we focus on methods for finding approximate inverse of $A$ by minimizing the Frobenius norm of the residual matrix $I-AP$. This minimization problem can be rewritten as \cite{Chow1998} %Let $P_1$ be good preconditioner for $A(1)$, i.e., computed by 
	\begin{align}\label{eq:gen-precond}
	\min_{P}\|I - AP\|_f^2.
	\end{align}
	Here, the columns of residual matrix $I-AP$ can be computed independently, which is an important property that can be exploited. Hence, the solution of (\ref{eq:gen-precond}) can be separated into $n$ independent least square problems as
	\begin{align}\label{eq:gen-precond-par}
	\begin{split}
	& \min_{P} \sum_{\mathcal{i}=1}^{n}\| (I-AP)e^{\mathcal{i}}\|_2^2, \ or \\
	& \min_{p^{\mathcal{i}}}\| e^{\mathcal{i}}-Ap^{\mathcal{i}}\|_2^2, \ \ \text{for} \ \mathcal{i} = 1, \ 2, \ \ldots, \ n,
	\end{split}
	\end{align} 
	where $e^{\mathcal{i}}$ and $p^{\mathcal{i}}$ are  the $\mathcal{i}$-th column of $I$ and  $P$, respectively. The above minimization problem can be implemented in parallel and one can efficiently obtain the explicit approximate inverse $P$ of $A$. 
	%	In Algorithm \ref{ALGO:AIRGASPAI}, we propose this SPAI preconditioner theory in the context of AIRGA algorithm. Here, only those parts that are relevant to solving the linear systems are listed. 	
}
%	Now, we apply this preconditioner theory (i.e., SPAI preconditioner) in AIRGA algorithm context, which is given in Algorithm~\ref{ALGO:AIRGASPAI}. Here, we only show those parts of AIRGA algorithm that require changes.  
%	\begin{algorithm}[!]
%		%	\large
%		%\normalsize
%					\setstretch{1.45}		
%		\caption{: AIRGA Algorithm with SPAI Preconditioner}
%		\label{ALGO:AIRGASPAI}
%		\begin{algorithmic}[1]
%			%\While {$||\hat{H}_{old}-\hat{H}_{new}||_{H_{2}} \ \leq \ tol$}
%			\State $z = 1$
%			\While{no convergence}
%			\For {$i = 1,\ \ldots,\ \ell$}  
%			\State Let $A_i^{(z)}=(s^{2}_{i}M+s_{i}D+K)$ 
%			\State Compute preconditioner $P_{i}^{(z)}$ by solving  
%			$\min\limits_{P_{i}^{(z)}}\|I-A_i^{(z)}P_{i}^{(z)}\|_{f}^{2}$
%			\State Solve $A_i^{(z)}P_{i}^{(z)}X^{(0)}(s_{i})=F$		 
%			\EndFor
%			%\While {$||\hat{H}_{j}-\hat{H}_{j+1}||_{H_{2}} \leq tol $}
%			\State j = 1
%			\While{no convergence and $j \le \lceil r_\text{max} / m \rceil $}
%			\For {$i = 1,\ \ldots, \ \ell$}
%			\State Only right hand sides are changing, so above preconditioner 
%			$P_{i}^{(z)}$ can \Statex \qquad \qquad \ \ be  applied as it is, i.e.,   
%			Solve $A_i^{(z)}P_{i}^{(z)}X^{(j)}(s_{i}) = -MV_{j}$
%			\EndFor
%			\State $j = j+1$
%			\EndWhile
%			\State $z = z+1$
%			\EndWhile			
%		\end{algorithmic}
%	\end{algorithm}

\subsection{Theory of reusing preconditioners}\label{sec:thr-resuse}
{\color{black}In general, the linear systems of equations generated by lines 4 and 10 of Algorithm \ref{ALGO:AIRGA}; lines 4 and 5 of Algorithm \ref{Algo:BIRKA}; and lines 4 and 10 of Algorithm \ref{Algo:QB-IHOMM} have the following form:
	\begin{align*}%\label{eq:gen-eqn}
	\begin{split}
	A_1X_1 & = F_1,\\
	A_2X_2 & = F_2,\\
	& \vdots \\
	A_{\ell}X_{\ell} & = F_{\ell},
	\end{split}
	\end{align*}     
	where $A_i \in \mathbb{R}^{n \times n}$, $X_i \in  \mathbb{R}^{n}$, and  $F_i \in  \mathbb{R}^{n}$; for $i = 1, \ 2,  \ \ldots, \ \ell$. 
	%	In our previous work \cite{navneet2019}, we presented two approaches of finding preconditioners. These two approaches are summarized in Table \ref{tab:cheap_update}. First approach works best for parametric dynamical systems, while second approach works well for non-parametric dynamical systems. In \cite{navneet2019}, we provided theoretical derivation and experiments for parametric dynamical systems only. Hence, we used first approach to obtain preconditioners. However, in this paper, we only focus on non-parametric dynamical systems. Therefore, we use second approach for finding a good preconditioners.   
	%	
	%	 As we have observed  in \cite{navneet2019}, there are two ways of finding preconditioners (see Table \ref{tab:cheap_update}), $1^{st}$ approach work best for parametric dynamical system, and $2^{nd}$ approach work well for non-parametric dynamical system. However, their we have provided theoretical derivation as well experiments is on parametric dynamical system only. In this paper our focus is on non-parametric dynamical system. Therefore, we apply $2^{nd}$ approach for finding a good preconditioners. 
	
	Let $P_1$ be a good preconditioner for ${A}_1$, that is, computed by  
	\begin{align*}
	\min_{P_1}\|I- {A}_1P_1\|_f^2.
	\end{align*}
	Now, we need to find a good preconditioner $P_2$ corresponding to ${A}_2$. {Using the standard SPAI theory, this means solving}
	\begin{align}\label{eq:first_precond}
	{\min_{P_2}\|I- {A}_2P_2\|_f^2.}
	\end{align}
	%However, there is an alternate cheaper way here. 
	If we are able to enforce ${A}_1P_1 = {A}_2 P_2$, then $P_2$ will be an equally good preconditioner for ${A}_2$ as much as $P_1$ is a good preconditioner for ${A}_1$ (since the Spectrum of ${A}_2P_2$ would be same as that of ${A}_1P_1$, on which convergence of any Krylov subspace method depends). Since $P_2$ is unknown here, we have a degree of freedom in choosing how to form it. Without loss of generality, we assume that $P_2 = Q_{2}P_1$, where $Q_{2}$ is an unknown matrix. Here, we need to enforce ${A}_1P_1 = {A}_2Q_{2}P_1$. {Thus, instead of solving the minimization problem (\ref{eq:first_precond}), we can solve}  
	\begin{align*}%\label{eq:second_precond}
	\min_{Q_2} \|{A}_1 - {A}_2 Q_{2}\|_f^2.
	\end{align*}
	%As we have mentioned \cite{navneet2019}, there are two ways of finding preconditioners (see Table \ref{tab:cheap_update}). In this table, $1^{st}$ approach work best for parametric case and $2^{nd}$ approach work well for non-parametric case. However, their we have provided theoretical derivation as well experiments is on parametric dynamical system only. As we have mentioned in this paper our focus is on non-parametric dynamical system. 
	%Hence, we apply a similar  
	%	Further, %we apply a similar
	%	%argument ($2^{nd}$ approach) for finding a good preconditioner 
	%	preconditioner for $\ell^{th}$  linear system, i.e., $P_{\ell}$ corresponding to ${A}_{\ell}$ is given below.%in the second column of  Table \ref{tab:cheap_update}.}
	%\begin{itemize}
	%	\item ${A}_{\ell-1}P_{\ell-1} = {A}_{\ell}P_{\ell}$
	%	\item If $P_{\ell} = Q_{\ell} P_{\ell-1} $,  then ${A}_{\ell-1}P_{\ell-1} = {A}_{\ell}Q_{\ell}P_{\ell-1}$
	%	\item $\min\limits_{Q_{\ell}} \|{A}_{\ell-1}-{A}_{\ell}Q_{\ell}\|_f^2$
	%\end{itemize}
	Note that $P_2$ here is never explicitly formed by multiplying two matrices $Q_2$ and $P_1$. Rather, always a matrix-vector product is done to apply the preconditioner. 
	
	Next, we apply a similar argument for finding a good preconditioner $P_i$ corresponding to $A_i$. For this we refer to one of our recent works \cite{navneet2019}, 
	which focused on MOR of parametric linear dynamical systems (category three from the Introduction).  
	We can obtain $P_i$ by enforcing either $A_{1}P_{1} = A_{i}P_i$ or $A_{i-1}P_{i-1} = A_iP_i$.  For these two cases, $P_i$ would be as effective preconditioner for $A_i$ as $P_1$ is for $A_1$ or $P_{i-1}$ is for $A_{i-1}$, respectively.  These two approaches are summarized  in Table \ref{tab:cheap_update}.
	\begin{table}[H]
		\centering
		\caption{{Cheap preconditioner update approaches \cite{navneet2019}.}}
		\setlength{\tabcolsep}{3pt}
		\def\arraystretch{1.0}
		\begin{tabular}{|l|l|}
			\hline
			First approach & Second approach \\ \hline
			%				\textbullet \ \ ${A}(1)P_1 = {A}(2)P_2$ &  \\
			%				\textbullet \ \ If $P_2 = Q_2P_1$, & \qquad Same as the \textit{first} approach\\
			%				\qquad then $ {A}(1)P_1 = {A}(2)Q_2P_1$ & \\
			%				\textbullet \ \ $\min\limits_{Q_2} \|{A}(1)-{A}(2) Q_2\|_f^2$ &   \\
			%				\hline
			%				\textbullet \ \ ${A}(1)P_1 = {A}(3)P_3$ & \textbullet \ \ ${A}(2)P_2 = {A}(3)P_3$ \\
			%				\textbullet \ \ If $P_3 = Q_3 P_1  $, & \textbullet \ \ If $P_3 = Q_3 P_{2}$, \\
			%				\qquad then $ {A}(1)P_1 = {A}(3)Q_3P_1$ & \qquad then ${A}(2)P_2 = {A}(3)Q_3P_2$ \\
			%				\textbullet \ \ $\min\limits_{Q_3} \| {A}(1)- {A}(3) Q_3\|_f^2$ & \textbullet \ \ $\min\limits_{Q_3} \|{A}{(2)}-{A}(3)Q_3\|_f^2$  \\
			%				\hline		
			%				\qquad \qquad \qquad \vdots & \qquad \qquad \qquad \vdots \\ \hline
			\textbullet \ \ ${A}_1P_1 = {A}_iP_i$ & \textbullet \ \ ${A}_{i-1}P_{i-1} = {A}_{i}P_{i}$ \\
			\textbullet \ \ If $P_{i} = Q_{i} P_1  $, & \textbullet \ \ If $P_{i} = Q_{i} P_{i-1} $, \\
			\qquad then $ {A}_1P_1 = {A}_{i}Q_{i}P_1$ & \qquad then ${A}_{i-1}P_{i-1} = {A}_{i}Q_{i}P_{i-1}$ \\
			\textbullet \ \ $\min\limits_{Q_{i}} \|{A}_1 - {A}_{i} Q_{i}\|_f^2$& \textbullet \ \ $\min\limits_{Q_{i}} \|{A}_{i-1}-{A}_{i}Q_{i}\|_f^2$  \\
			\hline	
		\end{tabular}
		\label{tab:cheap_update}
	\end{table}			
	
	In \cite{navneet2019}, we have conjectured (with evidence) the following two results:
	%\begin{itemize}
	(a) In the parametric case, the first approach is more beneficial. This is because, in this case although the two approaches  have a similarly hard minimization problem (attributed to slowly varying parameters, and in-turn, slowly changing matrices), the computation of $P_i$ from $P_1$ in the first approach leads to a preconditioner with less approximation errors, and hence, a one which is more accurate.
	(b) In the non-parametric case, the second approach is more suited. This is because in this case the minimization problem of the second approach is much easier to solve as compared to the first approach (attributed to rapidly changing expansion/ interpolation points, and in-turn, rapidly changing matrices). The computation of $P_i$ from $P_{i-1}$ in this case (rather than $P_1$ as above) does have the drawback of accumulated approximation errors, however, solving the minimization problem efficiently is a bigger bottleneck for scaling to large problems. 
	
	As mentioned in the Introduction, in \cite{navneet2019} we have extensively experimented for the parametric case (again, category three earlier) using the first approach. The focus here is to do a similar experimentation for the non-parametric case (first two categories earlier) using the second approach.
	%\end{itemize} 
}

\subsection{Application of reusing Preconditioner}\label{sec:app-reuse}
{\color{black} 
	Here, we first discuss the application of the above presented theory of reusing preconditioners to the AIRGA algorithm.  If we closely observe  Algorithm \ref{ALGO:AIRGA}, as mentioned earlier, linear systems are solved at lines 4 and 10. Computation of preconditioners is done only at line 4 because at line 10, matrices do not change, only the right-hand sides do. Hence, we only focus on reusing preconditioners for line 4. 
	
	Delving further into the complexity  of such linear systems, we observe that
	%, first time, the matrices corresponding to these linear systems changes with different expansion points (i.e., $s_i^2 M + s_i D + K$, where $i = 1, \ \ldots, \ \ell$). Second 
	the matrices 
	%corresponding to these linear systems again 
	change with the index of outer \texttt{while} loop (line 2) as well as with the index of the \texttt{for} loop corresponding to the expansion points (line 3). Hence, we denote such matrices not only with a subscript as in previous subsection but also with a superscript. That is, 
	$A_{i}^{(z)} = \left(s_i^{(z)}\right)^2M + s_i^{(z)} D + K, \ \text{where} \ z = 1, \ \ldots, \ \mathfrak{z} \ (\text{until covergence}) \ \text{and} \ i = 1, \ \ldots, \ \ell$. As the matrix $A_i^{(z)}$ changes with respect to two different indices, we can reuse preconditioners in many ways. However, here we use the second approach as discussed in the previous subsection. This approach is diagrammatically represented in Figure \ref{fig:reuse_AIRGA}. 
	%\vspace*{-.5cm}	
	\begin{figure}[h]
		\centering
		\includegraphics[width=95mm]{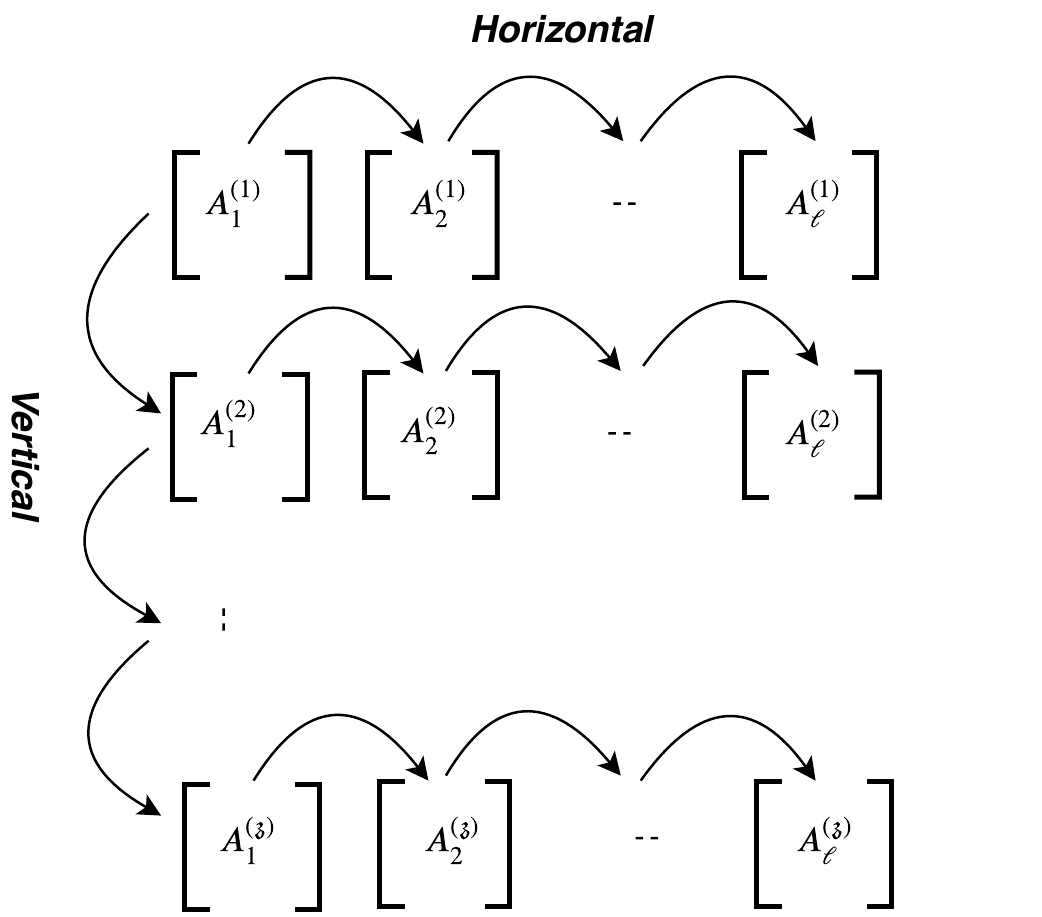}
		\caption{Reusing preconditioners in the AIRGA algorithm.}
		\label{fig:reuse_AIRGA}       % Give a unique label
	\end{figure}
	%\vspace*{-.24cm}
	%	Because of matrix changes at two indices, there are many ways of applying reuse of preconditioners, here, we apply the second approach as discussed above. The one which is best for us is in Figure \ref{fig:reuse_AIRGA}.
%	\Figure[!t](topskip=0pt, botskip=0pt, midskip=0pt){seq_matrix12.pdf}
%	{Reusing preconditioners in the AIRGA algorithm.\label{fig:reuse_AIRGA}}
	%\vspace*{-.4cm}
	%\vspace*{-.32cm}	
	%This approach works best for us because as discussed, reuse of preconditioners is more efficient when the change in the matrix is small. This phenomena is observed in the above applied technique. Hence, we propose to use this. 
	
	Next, we show how the new preconditioners are computed for both, the horizontal direction  and the vertical direction. While looking at the  horizontal route, 
	%discussed in Figure \ref{fig:reuse_AIRGA}. Next, we focus on horizontal direction. 	
	let, 	
	
	$A^{(z)}_{i-1} = \left(s_{i-1}^{(z)}\right)^2M + s_{i-1}^{(z)} D + K$ and 
	
	\ $A^{(z)}_{i} = \left(s_{i}^{(z)}\right)^2M + s_{i}^{(z)} D + K$ be the two coefficient matrices for different expansion points $s_{i-1}^{(z)}$ and $s_{i}^{(z)}$, respectively, with $i = 2, \ \ldots, \ \ell$. 
	%If we are able to enforce $A^{(1)}_1 P^{(1)}_1 = A^{(1)}_2 P^{(1)}_2$, then
	%	\begin{align*}
	%	A^{(1)}_1 P^{(1)}_1 \approx A^{(1)}_2 P^{(1)}_2.
	%	\end{align*}
	%	let $A_{{old}}=s_{old}^{2}M+s_{old}D+K$ and $A_{new}=s_{new}^{2}M+s_{new}D+K$ be the two coefficient matrices for different expansion points $s_{old}$ and $s_{new}$, respectively. These expansion points can be at the same or different AIRGA iterations. 
	%If the difference between $A_{old}$  and $A_{new}$ is small, then one can exploit this while building preconditioners for this sequence of matrices. 
	%This has been considered in the quantum Monte Carlo setting~\cite{ahuja2011improved} and for model reduction of first order linear dynamical systems~\cite{grim2015reusing,Wyatt2012}.  
	%	\par Let $P_{old}$ be a good initial  preconditioner for $A_{old}$. 
	%	As mentioned in earlier subsection, a reuse of preconditioner can be obtained  by enforcing $A_{old}P_{old} \approx A_{new}P_{new}$, where $old,\ new =\{1,\ \ldots, \ \ell\}$ and, as earlier, $\ell$ denotes the number of expansion points.    
	Using the above theory, we enforce $A^{(z)}_{i-1}P^{(z)}_{i-1} = A^{(z)}_{i} P^{(z)}_{i}$ in Figure \ref{fig:min_AIRGA}.
		\begin{figure*}[]
			\centering
			\includegraphics[]{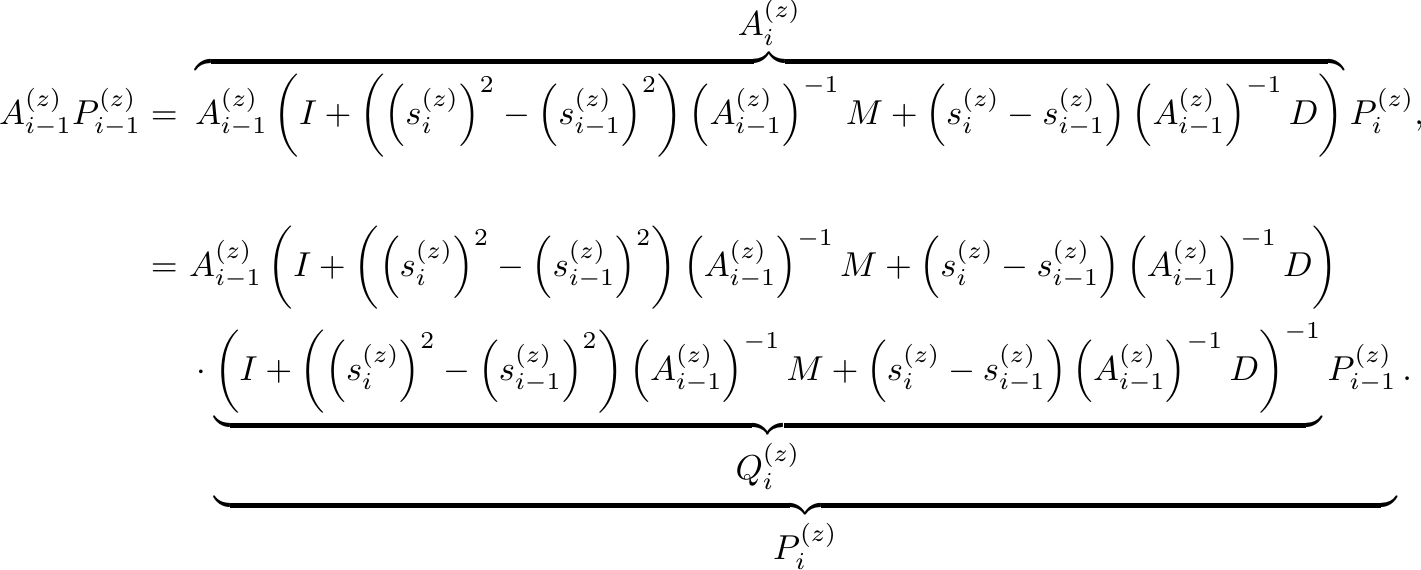}
			\caption{Expressing one linear system matrix in terms of the other.}
			\label{fig:min_AIRGA}       % Give a unique label
		\end{figure*}
	Thus, we eventually enforce $A^{(z)}_{i-1}P^{(z)}_{i-1} = A^{(z)}_{i} Q^{(z)}_{i} P^{(z)}_{i-1}$ and solve the  minimization problem
	\begin{equation*}
	\min_{Q^{(z)}_{i}}\|A^{(z)}_{i-1}-A^{(z)}_{i}Q^{(z)}_{i}\|^{2}_{f}. %= \min_{\left(q_2^1\right)^{(\mathcal{i})}} \sum_{\mathcal{i}=1}^{n}\|\left(a_1^1\right)^{(\mathcal{i})}-A_2^{(1)} \left(q_2^1\right)^{(\mathcal{i})}\|^{2}_2,
	\end{equation*}
%	where $\left(a_1^1\right)^{(\mathcal{i})}$ and $\left(q_2^1\right)^{(\mathcal{i})}$ denote the $\mathcal{i}^{th}$ columns of $A^{(1)}_1$ and $Q^{(1)}_2$, respectively.
This gives us the new preconditioner $P^{(z)}_i = Q^{(z)}_{i} P^{(z)}_{i-1}$.
This minimization is again performed for $n$ independent least square problems as in (\ref{eq:gen-precond-par}). 	
Similar steps are followed for reusing preconditioners along the rest of the horizontal directions, i.e. for all $z = 1, \ \ldots, \ \mathfrak{z}$.
	
	Now, applying this technique for the vertical direction, we  have $\text{for} \ z = 2, \ \ldots, \ \mathfrak{z}$	
	\begin{align*}
	A_1^{(z-1)}P_1^{(z-1)} = A_1^{(z)}P_1^{(z)}.
	\end{align*}	
	Following the steps as for the horizontal direction, here, we solve the minimization problem 
	\begin{align*}
	\min_{Q^{(z)}_{1}}\|A^{(z-1)}_{1}-A^{(z)}_{1}Q^{(z)}_{1}\|^{2}_{f}.
	\end{align*}
%	where
%	\begin{align*}
%	\resizebox{.48\textwidth}{!}
%	{$\begin{aligned}
%	Q^{(z)}_{1} = \left( I  + \left( \left(s_1^{(z)} \right)^2 - \left(s_1^{(z-1)} \right)^2  \right) M \left(A_1^{(z-1)}\right)^{-1} + 
%	%\right. \\	\left. 
%	\left( s_1^{(z)}  - s_1^{(z-1)}   \right) D \left(A_1^{(z-1)}\right)^{-1}\right)^{-1}.
%	\end{aligned}$
%    }
%	\end{align*} 
%	$Q^{(z)}_{1} = \left( I  + \left( \left(s_1^{(z)} \right)^2 - \left(s_1^{(z-1)} \right)^2  \right) M \left(A_1^{(z-1)}\right)^{-1} \right. \\ \left. \qquad \qquad \qquad \ \ \ +  
%	\left( s_1^{(z)}  - s_1^{(z-1)}   \right) D \left(A_1^{(z-1)}\right)^{-1}\right)^{-1}$. 
	This gives us the new preconditioner $P^{(z)}_1 = Q^{(z)}_{1} P^{(z-1)}_{1}$.
	Again, this is solved as $n$ independent least square problems as in (\ref{eq:gen-precond-par}).
	%	This is what we referred to as the second contribution, however, applying reusing preconditioners gets more challenging as the algorithm becomes intricate.  
	%	In Algorithm \ref{ALGO:AIRGA}, the preconditioner needs to be computed only at line 4 because no new preconditioner is computed at line 10. Hence, we only discuss the reuse of preconditioners for line 4. In  Algorithm \ref{ALGO:AIRGASPAI_Update}, we only give the code changes corresponding to lines 1 to 6 of Algorithm \ref{ALGO:AIRGA}. 
	%	%	Thus, application of this theory to the AIRGA given in  Algorithm \ref{ALGO:AIRGASPAI_Update}. In  Algorithm \ref{ALGO:AIRGASPAI_Update}, 
	%	 only the code corresponding line 1 to 6 changes. 
	
	AIRGA with an efficient implementation of the above discussed theory of reusing preconditioners is given in Algorithm \ref{ALGO:AIRGASPAI_Update}.	If we closely look at line 4 of Algorithm \ref{ALGO:AIRGA}, the solution vector is denoted by $X^{(0)}(s_i)$, where the superscript $``0"$ refers to the index of the inner \texttt{while} loop (line 8). We do not  bother about this index because, as earlier, matrix does not change inside this inner loop. Rather, we need to capture the change because of the outer \texttt{while} loop indexed with $z$.
	Hence, we denote the solution vector as $\mathcal{X}^{(z)}(s_i)$ in Algorithm \ref{ALGO:AIRGASPAI_Update} (lines 8, 11, 19 \& 22).
	%	of the reason discussed above. So, we simply denote it as $\mathcal{X}(s_i)$.
	%	 which we are does not bothered about at this stage because of the reason discussed above. Hence, we are denoting it as $\mathcal{X}(s_i)$.		
	It is important to emphasize again that preconditioners are never computed explicitly. 
	%It is always computed as products of previous preconditioners and the $Q$ matrices. However, this product is never explicitly formed by multiplying two matrices.
	Rather, they are obtained using matrix-vector product (please see line numbers 11, 19 \& 22 of Algorithm \ref{ALGO:AIRGASPAI_Update}). 
	%This efficient implementation is referred  to as our third contributions.  
	\begin{algorithm}[]
		%\large
		%		\normalsize
%		\setstretch{1.45}
		\caption{: AIRGA with reuse of SPAI preconditioner}
		\label{ALGO:AIRGASPAI_Update}
		{\color{black} \begin{algorithmic}[1]
				%\While {$||\hat{H}_{old}-\hat{H}_{new}||_{H_{2}} \ \leq \ tol$}
				\State $z = 1$
				\While{no convergence}
				%{\color{blue} \For {i = 1 to $l$}
				\If {$z == 1$}
				\For {$i = 1, \ \ldots, \  \ell$}
				\Statex
				\State $A_i^{(1)}= \left(\left(s_i^{(1)}\right)^{2}M + \left(s_i^{(1)}\right)D+K \right)$
				\Statex
				\If {$i == 1$}
				%\State $A_1^{(1)}=(s^{2}_1M+s_1D+K)$
				\State Compute initial $P_1^{(1)} $ by solving 
				\Statex \qquad \qquad \qquad \qquad $ \min\limits_{P_1^{(1)}}\|I-A_1^{(1)}P_1^{(1)}\|_{f}^2$ 
				\Statex \qquad \qquad \qquad {{\it \footnotesize{(First-time; no earlier preconditioner)}}}
				\Statex
				\State %Solve 
				$A_1^{(1)}P_1^{(1)}\mathcal{X}^{(1)}(s_1) = F$
				%\State Store $\mathcal{K}_{i_{old}} = \mathcal{K}_i$ and $P_{i_{old}} = P_i$
				\Else 
				\State Compute $Q_i^{(1)} $ by solving 
				\Statex \qquad \qquad \qquad \qquad $ \min\limits_{Q_i^{(1)}}\|A_{i-1}^{(1)}-A_i^{(1)}Q_i^{(1)}\|_{f}^2$ 
				\Statex \qquad \qquad \qquad  {\it \footnotesize{(Reuse along horizontal direction)}}
				\Statex
				\State %Solve 
				$A_{i}^{(1)} [{Q_{i}^{(1)} \ \cdots \ Q_2^{(1)}} P_1^{(1)}]\mathcal{X}^{(1)}(s_i) = F$
				\EndIf
				\EndFor
				\Else 
				\For {$i = 1, \ \ldots, \ \ell$}
				\Statex
				\State $A_i^{(z)} = \left(\left(s_i^{(z)}\right)^{2}M + \left(s_i^{(z)}\right)D+K \right)$
				\Statex
				\If {$i == 1$}
				%\State $A_1^{(1)}=(s^{2}_1M+s_1D+K)$
				\State Compute $Q_1^{(z)} $ by solving 
				\Statex \qquad \qquad \qquad \qquad $\min\limits_{Q_1^{(z)}}\|A_1^{(z-1)}-A_1^{(z)}Q_1^{(z)}\|_{f}^2$ 
				\Statex \qquad \qquad \qquad  {\it \footnotesize{(Reuse along vertical direction)}}
				\Statex
				%\State Solve 
				\State $A_1^{(z)}\left[Q_1^{(z)} \ \ldots \ Q_1^{(2)} \ P_1^{(1)} \right]\mathcal{X}^{(z)}(s_1) = F$
				%\State Store $\mathcal{K}_{i_{old}} = \mathcal{K}_i$ and $P_{i_{old}} = P_i$
				\Else 
				\State Compute $Q_i^{(z)} $ by solving 
				\Statex \qquad \qquad \qquad \qquad $ \min\limits_{Q_i^{(z)}}\|A_{i-1}^{(z)}-A_i^{(z)}Q_i^{(z)}\|_{f}^2$ 
				\Statex \qquad \qquad \qquad  {\it \footnotesize{(Reuse along horizontal direction)}} 
				\Statex
				%\State Solve 
				\State $A_{i}^{(z)} \left[\underbrace{Q_{i}^{(z)} \ \cdots \ Q_2^{(z)}} \right.$ \\ 
				\Statex \qquad \qquad \qquad \qquad $\left. \underbrace{Q_1^{(z)} \ \ldots \ Q_1^{(2)}} \  P_1^{(1)}\right] \mathcal{X}^{(z)}(s_i) = F$
				\EndIf
				\EndFor
				\EndIf
				\State  \text{{\it ``All the given set of expansion points}}
				\Statex \qquad \text{{(i.e. $s_1, \ s_2,  \ \ldots, \ s_{\ell}$) are updated''}}
				\State $z = z +1$
				\EndWhile	
				\Statex	{\it Note: The minimization problems at lines 7, 10, 18 and 21 are solved as $n$ independent least square problems (see (\ref{eq:gen-precond-par})).}
			\end{algorithmic}}
		\end{algorithm}
	} 
	% Compare this minimization problem with the one in SPAI (recall (\ref{eq:gen-precond})). Earlier, we were finding the preconditioner $P_{new}$ for $\mathcal{K}_{new}$ by solving $\min\limits_{P_{new}} \|I-\mathcal{K}_{new}P_{new}\|_{f}^{2}$. Here, we are finding the preconditioner $P_{new}$ (i.e., $P_{new}=Q_{new}P_{old}$) by solving $\min\limits_{Q_{new}}\|\mathcal{K}_{old}-\mathcal{K}_{new}Q_{new}\|_{f}^{2}$. The second formulation is much easier to solve, since in the first $\mathcal{K}_{new}$ could be very different from $I$, while in the second $\mathcal{K}_{new}$ and $\mathcal{K}_{old}$ are similar (change only in the expansion points). {\color{blue} For example, in our experiments section (Table \ref{tab:SPAi_SPAi_update_anly_4k} and Table \ref{tab:SPAi_SPAi_update_anly_12m}), we show that on an average the normed difference between $\mathcal{K}_{1}^{(z)}$ and $\mathcal{K}_{i}^{(z)}$ is many orders of magnitude smaller than the relative difference between $I$ and $\mathcal{K}_{i}^{(z)}$.}
	
	{\color{black} For sake of brevity, reusing preconditioners in BIRKA (Algorithm \ref{Algo:BIRKA}) is discussed as part of Appendix \ref{app:BIRKA}. Similarly, applying this theory to QB-IHOMM (Algorithm \ref{Algo:QB-IHOMM}) is discussed in Appendix \ref{app:QB-IHOMM}. }
	
%%%%%%%%%%%%%%%%%%%%%%%%%%%%%%%%%%%%%%%%%%%%%%%%%%%%%%%%%%%%%%
\section{Numerical Experiments}
\label{sec:Numerical_results}
{\color{black}
	For supporting our proposed preconditioned iterative solver theory using the AIRGA algorithm \cite{Bonin20161}, we perform experiments on two models. The first is a macroscopic equations of motion  model (i.e. academic disk brake $M_0$) \cite{Grabner2016}, and is discussed in Section \ref{sec:Acd_model}. The second is also a similar model, however, this is a real-life industrial problem (i.e. industrial disk brake $M_1$) \cite{Grabner2016}.
	%which is the $4^{th}$ contribution of this work as discussed in Section \ref{sec:intro}. 
	The experiments on this model are discussed in Section \ref{sec:Indu_model}.  
	%The second is macroscopic equations of motion (i.e., Industrial disk brake $M_1$) model \cite{Grabner2016}, and is discussed in Section \ref{sec:Indu_model}. 
	These models are described by the following set of equations \cite{Grabner2016}: 
	\begin{align}\label{eq:num}
	\begin{split}
	M_{\Omega} \ddot{x}(t)= & - D_{\Omega} \dot{x}(t) - K_{\Omega} x(t) + F u(t),  \\
	y(t) =  & \ C^T x(t),
	\end{split}  
	\end{align}  
	%where $M_{\Omega}, \ D_{\Omega}, \ K_{\Omega} \in \mathbb{R}^{n \times n}$ are parameter dependent coefficient matrices of mass, damping and stiffness, respectively. $F \in \mathbb{R}^n \ \text{and} \ C^T \in \mathbb{R}^n$ are taken as constant vectors (e.g., we have taken these vectors as a column of identity matrix).  Further, $M_{\Omega} = M, \ K_{\Omega} = K_E + K_R + \Omega^2 K_G, \ D_{\Omega} = \alpha M_{\Omega} + \beta K_{\Omega}$.  	
	%are parameter-dependent coefficient matrices, such as mass, damping and stiffness, respectively. 
	where $M_{\Omega} = M, \ K_{\Omega} = K_E + K_R + \Omega^2 K_G, \ D_{\Omega} = \alpha M_{\Omega} + \beta K_{\Omega}$ (case of proportionally damped system; as needed for AIRGA) with commonly used parameter values as ${\Omega} = 2 \pi, \ \alpha = 5 \times 10^{-02},$ and $\beta =  5 \times 10^{-06}$. Further, $F \in \mathbb{R}^n \ \text{and} \ C^T \in \mathbb{R}^n$ are taken as $\left[1 \ 0 \ \cdots \ 0 \right]^T$, which is the most frequently used choice. 
	%{\color{blue} We take four expansion points, i.e., $s_i = (1+3\frac{i-1}{4}),\ i = 1,\ldots, \ell$; $(\ell = 4)$ \cite{Grabner2016}.} 
	{\color{black} We take four expansion points linearly spaced between 1 and 500 based upon experience.}
	
	%	In these models, parameters are linearly embedded in the dynamical system matrices. However, the AIRGA algorithm is mainly proposed for non-parametric and especially for proportional damping. Therefore, we fix the given models according to our requirements, that is %$M_{\Omega} = M, \ K_{\Omega} = K_E + K_R + \Omega^2 K_G, D_{\Omega} = \alpha M_{\Omega} + \beta K_{\Omega}$, 
	%	by taking $\Omega = 2 \pi$, $\alpha = 0.05$ and $\beta = 0.000005$.
	%	%$F, \ C_p^T, \ \text{and} \ C_v^T \in \mathbb{R}^n$ are taken as constant vectors (e.g., we have taken these vectors as a column of identity matrix).
	%	%where $M,\ D,\ K \in  \mathbb{R}^{n \times n} $ are mass, damping and stiffness matrices, respectively, 
	%	%$F  \in  \mathbb{R}^{n \times m}, C_{p}, \\ C_{v}  \in  \mathbb{R}^{q \times n}$ are constant matrices.	
	%	%The range of variation of parameters $\Omega$ is [$\Omega_{min}$, $\Omega_{max}$] = $[1,4]*2\pi$. The coefficient matrices were constructed for a reference frequency $F_{ref} = 1600Hz$ and a reference angular velocity $\Omega = 1$ radians/sec \cite{Grabner2016}.
	%	We take four expansion points, i.e., $s_i = (1+3\frac{i-1}{4}),\ i = 1,\ldots, \ell$; $(\ell = 4)$ \cite{Grabner2016}. 
	
	Although our purpose is to just reuse SPAI in AIRGA (Algorithm \ref{ALGO:AIRGASPAI_Update}), we also execute original SPAI in AIRGA (Algorithm \ref{ALGO:AIRGA}) for comparison.	
	In Algorithms \ref{ALGO:AIRGA} and \ref{ALGO:AIRGASPAI_Update}, at line 2 the overall iteration (\texttt{while-loop}) terminates when the change in the reduced model (computed as $H_{2}$-error between the reduced models at two consecutive AIRGA iterations) is less than a certain tolerance. We take this tolerance as $10^{-04}$ based upon the values in~\cite{Bonin20161}. There is one more stopping criteria in Algorithms \ref{ALGO:AIRGA} at line 8 (also in Algorithm \ref{ALGO:AIRGASPAI_Update} but not listed here).  This checks the $H_{2}$-error between two temporary reduced models. We take this tolerance as $10^{-06}$, again based upon the values in~\cite{Bonin20161}. Since this is an adaptive algorithm, the optimal size of the reduced model is determined by the algorithm itself, and is denoted by $r$.

	The linear systems that arise here have non-symmetric matrices. 
	%As earlier, we use iterative methods instead of the direct method
	%\footnote{{\color{black}In our experiments, 
	%we have seen that the direct method is computationally expensive in the case of an industrial disk brake model (1.2 million in size ). A system 
	%		a machine with 16 GB RAM is unable to solve the linear systems arising from reducing the industrial brake model (1.2 million in size ) when
	%		using the direct method.  However, the preconditioned iterative method does work in this case.}} 
	%because the former scale well with increase in the size of the linear systems. 
	There are many iterative methods available for solving such linear systems. We use the Generalized Minimal Residual (GMRES) method \cite{Saad2003} because it is very popular \cite{Han2013}. The stopping tolerance in GMRES is taken as $10^{-06}$, which is a common standard. 
	%Preconditioning is employed when iterative methods fail or have slow convergence. Here, 
	As mentioned in Introduction, for both the given models, we observe that unpreconditioned GMRES fails to converge. Hence, we use the SPAI preconditioner as described above (without and with reuse). We use Modified Sparse Approximate Inverse (MSPAI 1.0) proposed in \cite{Alexander2008} as our preconditioner. This is because MSPAI uses a linear algebra library for solving sparse least square problems that arise here. We use standard initial settings of MSPAI $\left(\text{i.e.  tolerance (ep) of} \ 10^{-04} \right)$.% and cache size (cs) of 500.
}

{\color{black}We perform our numerical experiments on a machine with the following configuration: Intel Xeon (R) CPU E5-1620 V3 $@$ 3.50 GHz., frequency 1200 MHz., 8 CPU and 64 GB RAM. All the codes are written in MATLAB (2016b) (including AIRGA, GMRES) except SPAI and reusable SPAI. MATLAB is used because of ease of rapid prototyping. Computing SPAI and reusable SPAI in MATLAB is expensive, therefore, we use C++ version of these (SPAI is from MSPAI and reusable SPAI is written by us). MSPAI further uses BLAS, LAPACK and ATLAS libraries. 
	%While solving every linear system in the sequence, 
	Whenever a preconditioner has to be computed, we first compute SPAI and reusable SPAI  separately (in-parallel) and save them. Then, we run MATLAB code along with the saved preconditioner matrices (i.e. SPAI and reusable SPAI).}

%%%%%%%%%%%%%%%%%%%%%%%%%%%%%%%%%%%%%%%%%%%%%%%%%%%%%%%
\subsection{Academic Disk Brake Model} \label{sec:Acd_model} 
{\color{black}This model is of size $4,669$. Based upon experience, the maximum reduced system size ($r_{max}$) is taken as $20$. As mentioned earlier, however, due to the adaptive nature of the AIRGA algorithm, we obtain a reduced system of size $r = 13$. For this model, the AIRGA algorithm takes two outer iterations (line 2 of Algorithms \ref{ALGO:AIRGA} and \ref{ALGO:AIRGASPAI_Update}) to converge (i.e. $\mathfrak{z}=2$).
	%We compute a reduced order model  by the AIRGA algorithm~\cite{Bonin20161}.
	%Here, $r_{max}$ is taken as $20$ (reduced system size). However, we obtained the reduced  system with size $r = 13$ due to adaptive nature of AIRGA algorithm.
	% We take $r_\text{max}$, i.e., the maximum dimension to which we want $13$. 

	Reusing the SPAI preconditioner is beneficial  %Recall, the first linear system matrix defined as
	%	\begin{align*}
	%	A_1^{(1)} = \left(s_1^{(1)}\right)^2 M + \left(s_1^{(1)}\right) D + K, 
	%	\end{align*}
	%	and other matrices,
	%	\begin{align*}
	%	A_i^{(z)} = \left(s_i^{(z)}\right)^2 M + \left(s_i^{(z)}\right) D + K, 
	%	\end{align*}
	%	for $i = 1, \ldots, \ell$ with $\ell = 4$. Based upon our proposed theory of reusing preconditioner in Section \ref{sec:app-reuse} (see Algorithm \ref{ALGO:AIRGASPAI_Update}, Table \ref{tab:cheap_update} and Figure \ref{fig:reuse_AIRGA}), reusing SPAI is useful 	
	when the values of $\|I-A_i^{(z)}\|_f/\|I\|_f$ is large, and the values of $\|A_{i-1}^{(z)}-A_i^{(z)}\|_f/\|A_{i-1}^{(z)}\|_f$ and $\|A_1^{(z-1)}-A_1^{(z)}\|_f/\|A_1^{(z-1)}\|_f$ are small,  which is true in this case (see Table \ref{tab:SPAi_SPAi_update_anly_4k}). In this table, columns 1 and 2 list the AIRGA iterations and the four expansion points, respectively. The above three quantities are listed in columns 3, 4 and 5, respectively.     
	%This can be observed from Table \ref{tab:SPAi_SPAi_update_anly_4k} (see columns 3, 4 and 5). In this table, we have given both the cases (horizontal and vertical directions). 
	For the first AIRGA iteration and the first expansion point, SPAI preconditioner cannot be reused because there is no earlier preconditioner (mentioned as {\it NA} in table). 
	From the second expansion point (and the first AIRGA iteration), we perform horizontal reuse of preconditioner (see Figure \ref{fig:reuse_AIRGA}). This is the same for the second AIRGA iteration as well. Vertical reuse of preconditioner is done only for the first expansion point (and the second AIRGA iteration; again see Figure \ref{fig:reuse_AIRGA}). 
	%These details are given in columns 4 and 5. 
	
	%This is the case of horizontally  reused of SPAI preconditioner. Next, we check how close these matrices are in vertically. For that we only check first linear system matrices in each of the AIRGA iterations (here we have two AIRGA iterations). That is, $\|A_1^{(1)}-A_1^{(2)}\|_f/\|A_1^{(1)}\|_f = 0.9996.$ This is also small as compared to  $\|I-A_1^{(2)}\|_f/\|I\|_f$. Thus, we use both the strategies (i.e., horizontal and vertical) together to get the benefits of reused SPAI.      
	
	%As discussed in Section~\ref{subsec:SPAI_Update}, SPAI update in AIRGA is done horizontally (see Algorithm~\ref{ALGO:AIRGASPAI_Update} and Table~\ref{change_in_expansion_pnt}). Therefore, here, we are computing preconditioner at each AIRGA iteration with the first expansion point only. Preconditioners for the remaining expansion points in first AIRGA iteration, and other AIRGA iterations are  computed using SPAI update (i.e., $||\mathcal{K}^{(z)}_1-\mathcal{K}^{(z)}_i||_f$).
	
	In Table~\ref{tab:SPAI_SPAI_Update_time_4k}, we compare the SPAI and the reusable SPAI timings.  As for Table \ref{tab:SPAi_SPAi_update_anly_4k}, here columns 1 and 2 list the AIRGA iterations and the four expansion points, respectively. SPAI and reusable SPAI computation times are given in columns 3 and 4, respectively.  
	%	In the first column, we give the AIRGA iterations. The number of expansion points in each AIRGA iterations are given in second column. The computation time of SPAI and reusing SPAI are give in columns 3 and 4, respectively. 
	At the first AIRGA iteration and the first expansion point, both SPAI and reusable SPAI take the same computation time.  This is because, as above, reusing of SPAI preconditioner is not applicable here. 	%from second expansion points onwards of the first AIRGA iteration. %, saving in time is observed from these steps.
	%This is because SPAI update is applicable from second expansion points of first AIRGA iterations onwards and also second expansion points of second AIRGA iterations onwards, saving in time is observed from these step onwards. 
	From the second expansion point of the first AIRGA iteration,
	% we can be observed that, while using reused SPAI, it reduces the preconditioner computation time from $22$ mins to $7$ mins (approximately). This shows a considerable 
	we see substantial savings because of the reuse of the SPAI preconditioner (approximately $68 \%$). %As discussed in earlier paragraphs, %since SPAI update is being done only from AIRGA iteration 3, 	
	\begin{table}[]
		%\hspace{-1cm}
		%\begin{minipage}[b]{.5\textwidth}
		%\centering	
		\def\arraystretch{1.4}
		\setlength\tabcolsep{1.5pt}
		\caption{{SPAI and reusable SPAI analysis for the academic disk brake model.}}
		\label{tab:SPAi_SPAi_update_anly_4k}
	\begin{threeparttable}
		\begin{tabular}{|c|c|c|c|c|}
			\hline
			\multirow{4}{*}{\begin{tabular}[c]{@{}c@{}} AIRGA  \\ Itr.$^\dagger$\end{tabular}} & \multirow{4}{*}{\begin{tabular}[c]{@{}c@{}} {Exp.} \\ {Pts.}$^\ddagger$ \end{tabular}} & \begin{tabular}[c]{@{}c@{}} {SPAI Case} \end{tabular} & \multicolumn{2}{c|}{\begin{tabular}[c]{@{}c@{}} {Reusable SPAI Case} \end{tabular}} \\ \cline{3-3} \cline{4-4} \cline{5-5}			
			\multirow{8}{*}{1} &  & {\begin{tabular}[c]{@{}c@{}} {Standard} \vspace{0.2cm} \\ $\footnotesize{\dfrac{\|I - A_i^{(z)}\|_f}{\|I\|_f}}$ \end{tabular}} &  {\begin{tabular}[c]{@{}c@{}} {Horizontal} \vspace{0.2cm} \\ $\footnotesize{\dfrac{\| A_{i-1}^{(z)} - A_i^{(z)} \|_f}{\| A_{i-1}^{(z)}  \|_f}}$ \end{tabular}}  & {\begin{tabular}[c]{@{}c@{}} {Vertical} \vspace{0.2cm} \\ $ \footnotesize{\dfrac{\| A_1^{(z-1)} - A_1^{(z)} \|_f}{\| A_1^{(z-1)}  \|_f}}$ \end{tabular}}  \\ \cline{1-1} \cline{2-2} \cline{3-3} \cline{4-4} \cline{5-5} 
			& 1 &  $3.77 \times 10^{06}$ & {\it NA} & {\it NA} \\ \cline{2-2} \cline{3-3} \cline{4-4} \cline{5-5}
			& 2 &  $4.36 \times 10^{06}$ & $0.1569$ &  \\ \cline{2-2} \cline{3-3} \cline{4-4} 
			& 3 &  $4.95 \times 10^{06}$ & $0.3139$ & {\it NA} \\ \cline{2-2} \cline{3-3} \cline{4-4} 
			& 4 &  $5.54 \times 10^{06}$ & $0.4708$ &  \\ \cline{1-1} \cline{2-2} \cline{3-3} \cline{4-4} \hline
			\multirow{4}{*}{2} & 1 & $7.63 \times 10^{06}$ & {\it NA}  & $0.9996$ \\ \cline{2-2} \cline{3-3} \cline{4-4} \cline{5-5}
			& 2 &  $4.06 \times 10^{06}$ &  $0.0180$ &  \\ \cline{2-2} \cline{3-3} \cline{4-4} 
			& 3 &  $1.62 \times 10^{06}$ &  $20.3431$ & {\it NA}  \\  \cline{2-2} \cline{3-3} \cline{4-4} 
			& 4 &  $3.82 \times 10^{06}$ &  $0.4985$ &   \\ \hline 
		\end{tabular}
		\begin{tablenotes}
			\footnotesize{\item[$\dagger$] AIRGA Iterations. \item[$\ddagger$] Expansion Points}.
		\end{tablenotes}
	\end{threeparttable}		
	\end{table}
	%\end{minipage}\qquad \qquad \qquad \qquad 
	\begin{table}[]
		%\begin{minipage}[b]{.5\textwidth}
		\centering
		\def\arraystretch{0.96}
		\setlength\tabcolsep{1.5pt}
		\caption{ SPAI and reusable SPAI computation time for the academic disk brake model.}
		\label{tab:SPAI_SPAI_Update_time_4k}
		\begin{tabular}{|c|c|c|c|}%{p{3.5 cm} p{8 cm} p{5 cm}}
			\hline
			{{\begin{tabular}[c]{@{}c@{}} AIRGA \\Iterations $(z)$ \end{tabular}}} & {\begin{tabular}[c]{@{}c@{}} Expansion \\ Points $(s_i)$\end{tabular}}&{ \begin{tabular}[c]{@{}c@{}} SPAI \\ (Seconds)  \end{tabular}} & {\begin{tabular}[c]{@{}c@{}} Reusable SPAI \\ (Seconds) \end{tabular}} \\ \hline
			\multirow{4}{*}{1} & 1 & 174 & 174 \\ \cline{2-2} \cline{3-3} \cline{4-4} 
			& 2 & 164  & 10 \\ \cline{2-2} \cline{4-4} \cline{3-3} 
			& 3 & 165 &  16 \\ \cline{2-2} \cline{4-4} \cline{3-3} 
			& 4 & 165 &  20 \\ \cline{1-1}\cline{2-2} \cline{3-3}  \cline{4-4} \hline
			\multirow{4}{*}{2}  & 1 & 165  & 64 \\ \cline{2-2} \cline{3-3} \cline{4-4} 
			&  2 & 165 &  10 \\  \cline{2-2}  \cline{3-3} \cline{4-4} 
			&  3 & 165 &  108 \\  \cline{2-2}   \cline{3-3} \cline{4-4} 
			&  4 & 158 &  20 \\ \hline 
			\cline{1-4} \bigstrut
			\ \ \ \ \ \textbf{Total}	& {\bf 8} & {\bf 1321} &  {\bf 422} \\ \hline		
		\end{tabular}		
		%\end{minipage}
	\end{table}
		
	Table~\ref{tab:GMRES4K} provides the iteration count and the computation time of GMRES. Here, we only provide GMRES execution details since the computation time of preconditioner has been discussed above. 
	%\footnote{\color{blue}{Here, the comparison between GMRES with SPAI and  GMRES with reusing SPAI are not given. The reason is, in both cases, the differences are very less, GMRES with reusing SPAI takes slightly more time (i.e., in a fraction of second) as compare to GMRES with SPAI.}}.
	%	In the first column, again we give AIRGA iterations. The number of times GMRES is executed for each of the linear systems are given in the second column. 
	%	The iteration count and computation time of GMRES are given in columns 3 and 4, respectively.
	In this table, column 1 lists the AIRGA iterations. The number of linear solves and average GMRES iterations per linear solve are given in columns 2 and 3, respectively. Finally, columns 4 and 5 list the computation times of GMRES when using SPAI and reusable SPAI, respectively.   
	We notice from this table that solving linear systems by GMRES with SPAI takes less computation time as compared to solving them by GMRES with reusable SPAI. This is because when we reuse the SPAI preconditioner in GMRES, additional matrix-vector products are performed, however, this extra cost is almost negligible when compared to the savings in the preconditioner computation time for the latter case (as evident in Table \ref{tab:SPAi_SPAi_update_anly_4k} above; also see total GMRES and preconditioner time below).  
	%Here, we have not added the computation time of preconditioner. As we have already mentioned earlier, unpreconditioned GMRES does not converge for the given model. %Here, we have taken maximum iteration in GMRES for convergence is $500$.% and $10^{-06}$ as a stopping tolerance. 
	
	Table~\ref{tab:gmres_spai_spai_u4k} gives the computation time of GMRES plus SPAI (column 2) and GMRES plus reusable SPAI (column 3) at each AIRGA iteration (column 1). As evident from this table, reusing the SPAI preconditioner leads to about $60 \%$ savings in total time required for solving all the linear systems. 
	%	As earlier, the AIRGA iterations are given in first column. 
	%	These admissible computation times are given in columns 2 and 3, respectively. 
	%	% when using basic SPAI preconditioner  and when using SPAI with update preconditioner  for the given model. 
	%	%We can notice from this table that iterative solves with update take about ? of time as needed for basic iterative solves (or without SPAI update).
	%	We can observe from this table that computation time of GMRES with reusing SPAI is on an average $58 \%$ less than the computation time of GMRES with SPAI.
} 
\begin{table*}[]
	%	\begin{minipage}[b]{.44\textwidth}
	\centering
			\def\arraystretch{1.0}
		\setlength\tabcolsep{3pt}
	\caption{GMRES computation time for the academic disk brake model.}
	\label{tab:GMRES4K}
	\begin{tabular}{|c|c|c|c|c|}
		\hline
		\multirow{2}{*}{{\begin{tabular}[c]{@{}c@{}}AIRGA \\ Iterations $(z)$\end{tabular}}}   
		&{\begin{tabular}[c]{@{}c@{}} No. of \\ Linear Solves \end{tabular}}& {\begin{tabular}[c]{@{}c@{}} GMRES Iterations \\ per Linear Solve \end{tabular}}& {\begin{tabular}[c]{@{}c@{}} GMRES Time when \\ Using SPAI \\  (Seconds)\end{tabular}} & {\begin{tabular}[c]{@{}c@{}} GMRES Time when \\ Using Reusable SPAI \\  (Seconds)\end{tabular}} \\ \hline
		1 & 10 & 271 & $7.98 $  & $8.62$ \\ \hline
		2 & 13 & 270 & $8.12 $  & $8.93$ \\ \hline
		\cline{1-5} \bigstrut[b]
		\ \ \ \ \ \textbf{Total}& $\mathbf{23}$ & {\bf \begin{tabular}[c]{@{}c@{}}  ${10 \times 271 + 13 \times 270}$ 
				\\ $\mathbf{= 6220}$ \end{tabular}} &\textbf{\begin{tabular}[c]{@{}c@{}} ${10 \times 7.98 + 13 \times 8.12}$\\ $\mathbf{=  185}$ \end{tabular}} & \textbf{\begin{tabular}[c]{@{}c@{}} ${10 \times 8.62 + 13 \times 8.93}$\\ $\mathbf{=  202}$ \end{tabular}} \\ \hline
	\end{tabular}
\end{table*}
%	\end{minipage}\qquad \qquad
%	\begin{minipage}[b]{.44\textwidth}
\begin{table}[]
	\centering
\setlength\tabcolsep{3pt}
			\def\arraystretch{1.0}
	\caption{GMRES with SPAI and reusable SPAI computation time for the academic disk brake model.}
	\label{tab:gmres_spai_spai_u4k}
	\begin{tabular}{|c|c|c|}
		\hline
		{\begin{tabular}[c]{@{}c@{}}AIRGA \\ Iterations $(z)$\end{tabular}} & {\begin{tabular}[c]{@{}c@{}}GMRES Plus\\ SPAI Time \\ (Seconds)\end{tabular}} & {\begin{tabular}[c]{@{}c@{}}GMRES Plus \\ Reusable SPAI Time \\ (Seconds) \end{tabular}} \\ \hline
		1 & 748  & 306 \\ \hline
		2 & 759  & 318    \\ \hline
		\cline{1-3} \bigstrut
		\ \ \ \ \ \textbf{Total} & {\bf 1507} & {\bf 624}  \\ \hline
	\end{tabular}
	%	\end{minipage}
\end{table}

%%%%%%%%%%%%%%%%%%%%%%%%%%%%%%%%%%%%%%%%%%
\subsection{Industrial Disk Brake Model}\label{sec:Indu_model}
{\color{black}%This model is a industrial brake model of size 1.2 million (i.e., $M_1$ model of \cite{Grabner2016}). 
	%	 As we have mentioned in Section \ref{sec:intro}, there are four contributions to this work. Here, we discuss our $4^{th}$ contribution. We experiment on real-life industrial problem of size 1.2 million (i.e., industrial disk brake model; $M_1$ model of \cite{Grabner2016} ). 	
	%	%We compute a reduced order model  by the AIRGA algorithm~\cite{Bonin20161}. %We take $r_\text{max}$, i.e., the maximum dimension to which we want $300$ \cite{Grabner2016}. 
	%	Here, $r_{max}$ is taken as $100$ (reduced system size). However, we obtained the reduced  system with size $r = 52$ because of the reason mentioned earlier.
	
	This model is of size $1.2$ million. Based upon experience, the maximum reduced system size ($r_{max}$) is taken as $100$. As mentioned earlier, however, due to the adaptive nature of the AIRGA algorithm, we obtain a reduced system of size $r = 52$. For this model, the AIRGA algorithm takes four outer iterations (line 2 of Algorithms \ref{ALGO:AIRGA} and \ref{ALGO:AIRGASPAI_Update}) to converge (i.e. $\mathfrak{z}=4$).

	Again, reusing the SPAI preconditioner is beneficial when the value of $\|I-A_i^{(z)}\|_f/\|I\|_f$ is large, and the value of $\|A_{i-1}^{(z)}-A_i^{(z)}\|_f/\|A_{i-1}^{(z)}\|_f$ and $\|A_1^{(z-1)}-A_1^{(z)}\|_f/\|A_1^{(z-1)}\|_f$ are small,  which is true in this case (see Table \ref{tab:SPAi_SPAi_update_anly_12m}). The structure of this table is same as Table \ref{tab:SPAi_SPAi_update_anly_4k}.	As earlier, for the first AIRGA iteration and the first expansion point, SPAI preconditioner cannot be reused because there is no earlier preconditioner (mentioned as {\it NA} in table). From the second expansion point (and the first AIRGA iteration), we perform horizontal reuse of preconditioner (see Figure \ref{fig:reuse_AIRGA}). This is the same for the second, the third and the fourth AIRGA iterations as well. Vertical reuse of preconditioner is done only for the first expansion point (and the second, the third, and the fourth AIRGA iterations; again see Figure \ref{fig:reuse_AIRGA}).  
	\begin{table*}%
		\centering
		\setlength\tabcolsep{1.5pt}
		\parbox{0.5\textwidth}{
			\caption{{SPAI and reusable SPAI analysis for the industrial disk brake model.}}
			\begin{small}
				\def\arraystretch{1.4}
				\begin{tabular}{|c|c|c|c|c|}
					\hline
					\multirow{4}{*}{\begin{tabular}[c]{@{}c@{}} AIRGA  \\ Iterations $(z)$ \end{tabular}} & \multirow{4}{*}{\begin{tabular}[c]{@{}c@{}} {Expansion} \\ {Points}  $(s_i)$ \end{tabular}} & {SPAI Case} & \multicolumn{2}{c|}{\begin{tabular}[c]{@{}c@{}} {Reusable SPAI Case} \end{tabular}} \\ \cline{3-3} \cline{4-4} \cline{5-5}			
					\multirow{8}{*}{1} &  & {\begin{tabular}[c]{@{}c@{}} {Standard} \vspace{0.2cm} \\  $\footnotesize{\dfrac{\|I - A_i^{(z)}\|_f}{\|I\|_f}}$ \end{tabular}} &  {\begin{tabular}[c]{@{}c@{}} {Horizontal} \vspace{0.2cm} \\  $\footnotesize{\dfrac{\| A_{i-1}^{(z)} - A_i^{(z)} \|_f}{\| A_{i-1}^{(z)}  \|_f}}$ \end{tabular}}  & {\begin{tabular}[c]{@{}c@{}} {Vertical} \vspace{0.2cm} \\  $ \footnotesize{\dfrac{\| A_1^{(z-1)} - A_1^{(z)} \|_f}{\| A_1^{(z-1)}  \|_f}}$ \end{tabular}}  \\ \cline{1-1} \cline{2-2} \cline{3-3} \cline{4-4} \cline{5-5} 
					& 1 &  $6.54 \times 10^{08}$  & {\it NA} & {\it NA} \\ \cline{2-2} \cline{3-3} \cline{4-4} \cline{5-5}
					& 2 &  $6.54 \times 10^{08}$  & $3.74 \times 10^{-05}$ & \\ \cline{2-2} \cline{3-3} \cline{4-4} 
					& 3 &  $6.54 \times 10^{08}$  & $7.49 \times 10^{-05}$ & {\it NA} \\ \cline{2-2} \cline{3-3} \cline{4-4} 
					& 4 &  $6.54 \times 10^{08}$  & $1.12 \times 10^{-04}$ & \\ \cline{1-1} \cline{2-2} \cline{3-3} \cline{4-4} \hline
					\multirow{4}{*}{2} & 1 & $1.31 \times 10^{09}$  & {\it NA} & $1.006$ \\ \cline{2-2} \cline{3-3} \cline{4-4} \cline{5-5}
					& 2 &  $6.65 \times 10^{08}$ &  $0.49$ &\\ \cline{2-2} \cline{3-3} \cline{4-4} 
					& 3 &  $6.53 \times 10^{08}$ &  $0.50$ & {\it NA} \\  \cline{2-2} \cline{3-3} \cline{4-4} 
					& 4 &  $6.56 \times 10^{08}$ &  $0.49$  &\\ \hline 
					\multirow{4}{*}{3} %&  & $\|I - \mathcal{K}_i^{(z)}\|_f/\|I\|_f$   & $\| \mathcal{K}_1^{(1)} - \mathcal{K}_i^{(z)} \|_f/\| \mathcal{K}_1^{(1)}  \|_f$ \\ \cline{2-2} \cline{3-3} \cline{4-4} 
					& 1 &  $1.30 \times 10^{09}$  & {\it NA}  & $1.009$ \\ \cline{2-2} \cline{3-3} \cline{4-4} \cline{5-5}
					& 2 &  $7.01 \times 10^{08}$  & $0.4658$ &\\ \cline{2-2} \cline{3-3} \cline{4-4} 
					& 3 &  $6.53 \times 10^{08}$  & $0.5499$ & {\it NA}\\ \cline{2-2} \cline{3-3} \cline{4-4} 
					& 4 &  $6.63 \times 10^{08}$  & $0.4940$ &\\ \cline{1-1} \cline{2-2} \cline{3-3} \cline{4-4} \cline{5-5} \hline
					\multirow{4}{*}{4} & 1 & $1.31 \times 10^{09}$  & {\it NA} & $1.0015$\\ \cline{2-2} \cline{3-3} \cline{4-4}  \cline{5-5}
					& 2 &  $6.86 \times 10^{08}$ &  $0.4641$ &\\ \cline{2-2} \cline{3-3} \cline{4-4} 
					& 3 &  $6.53 \times 10^{08}$ &  $0.5002$ & {\it NA} \\  \cline{2-2} \cline{3-3} \cline{4-4} 
					& 4 &  $6.56 \times 10^{08}$ &  $0.4933$ & \\ \hline 
				\end{tabular}
			\end{small}
			\label{tab:SPAi_SPAi_update_anly_12m}}
		\centering
		\qquad \qquad
			\parbox{0.3\textwidth}{
				\caption{{SPAI and reusable SPAI  computation time for the industrial disk brake model.}}
		\begin{small}
			\def\arraystretch{1.4}
			\begin{threeparttable}
			\begin{tabular}{|c|c|c|c|}%{p{3.5 cm} p{8 cm} p{5 cm}}
				\hline
				{\begin{tabular}[c]{@{}c@{}} AIRGA \\ Iterations $(z)$ \end{tabular}}& {\begin{tabular}[c]{@{}c@{}} Expansion \\ Points $(s_i)$\end{tabular}}& {\begin{tabular}[c]{@{}c@{}} SPAI\tnote{$\S$}  \end{tabular}} & {\begin{tabular}[c]{@{}c@{}} Reusing SPAI\tnote{$\S$}  \end{tabular}} \\ \hline
				\multirow{4}{*}{1} & 1 & 10 \text{hrs}  & 10 \text{hrs}  \\ \cline{2-2} \cline{3-3} \cline{4-4} 
				& 2 & 10 \text{hrs}   & 1 \text{hr}  \\ \cline{2-2} \cline{4-4} \cline{3-3} 
				& 3 & 10 \text{hrs}  & 1 \text{hr} \\ \cline{2-2} \cline{4-4} \cline{3-3} 
				& 4 & 10 \text{hrs}  & 1 \text{hr} \\ \cline{1-1}\cline{2-2} \cline{3-3}  \cline{4-4} \hline
				\multirow{4}{*}{2}  & 1 & 10 \text{hrs}  & 1 \text{hr} 30 \text{mins} \\ \cline{2-2} \cline{3-3} \cline{4-4} 
				&  2 & 10 \text{hrs}  &  1 \text{hr} \\  \cline{2-2}  \cline{3-3} \cline{4-4} 
				&  3 & 10 \text{hrs}  &  1 \text{hr} \\  \cline{2-2}   \cline{3-3} \cline{4-4} 
				&  4 & 10 \text{hrs}  & 1 \text{hr} \\ \hline 
				\multirow{4}{*}{3} & 1 & 10 \text{hrs}  & 1 \text{hr} 30 \text{mins} \\  \cline{2-2} \cline{3-3} \cline{4-4} 
				& 2 & 10 \text{hrs}  & 1  \text{hour} \\ \cline{2-2} \cline{3-3} \cline{4-4} 
				& 3 & 10 \text{hrs}   & 1 \text{hour} \\  \cline{2-2} \cline{3-3} \cline{4-4} 
				& 4 & 10 \text{hrs}  & 1  \text{hour} \\ \cline{1-1}\cline{2-2} \cline{3-3} \cline{4-4} 
				\multirow{4}{*}{4} & 1 & 10 \text{hrs} & 1 hr 30 \text{mins} \\ \cline{2-2} \cline{3-3} \cline{4-4} 
				& 2 & 10 \text{hrs}  & 1 \text{hr} \\ \cline{2-2} \cline{3-3} \cline{4-4} 
				& 3 & 10 \text{hrs}  & 1 \text{hr} \\  \cline{2-2} \cline{3-3} \cline{4-4} 
				& 4 & 10 \text{hrs}  & 1 \text{hr} \\ \hline
				\cline{1-4} \bigstrut
				\ \ \ \ \ \textbf{Total}	& {\bf 16} & {\bf 160 \text{hrs}} & {\bf 26 \text{hrs} 30 \text{mins}} \\ \hline		
			\end{tabular}
			\begin{tablenotes}
				\footnotesize{\item[$\S$] All times given here differ in seconds (not evident because of 
					the rounding to the nearest minute)}.
			\end{tablenotes}
		\end{threeparttable}
		\end{small}
		\label{tab:SPAI_SPAI_Update_time_12m}}
	\end{table*}

	In Table~\ref{tab:SPAI_SPAI_Update_time_12m}, we compare the SPAI and the reusable SPAI timings. The structure of this table is same as that of Table \ref{tab:SPAI_SPAI_Update_time_4k}.
	%	In the first column, we give the AIRGA iterations. The number of expansion points in each AIRGA iteration is given in the second column. The computation time of SPAI and reusing SPAI are given in columns 3 and 4, respectively. 
	As before, at the first AIRGA iteration and the first expansion point, both SPAI and reusable SPAI take the same computation time.  This is because, as above, reusing of SPAI preconditioner is not applicable here. From the second expansion point of the first AIRGA iteration, we see substantial savings because of the reuse of the SPAI preconditioner (from $160$ hours to $26$ hrs $30$ minutes; approximately $83 \%$).
	
	%	Again, at the first AIRGA iteration and first expansion point, both SPAI and reusing SPAI takes the same computation time. This is because reusing SPAI is not applicable at this stage. From the second expansion point of first AIRGA iteration, we see substantial savings in computation time (from $160$ hrs to $26$ hrs $30$ mins; approximately $83 \%$) with reusing SPAI as compared to SPAI. %This shows the considerable amount of time is saved  by using SPAI with updates (i.e., $83 \%$). %As discussed in earlier paragraphs, %since SPAI update is being done only from AIRGA iteration 3, saving in time is observed from this step onwards only. 
	Table~\ref{tab:GMRES} provides the iteration count and the computation time of GMRES. Here, again we have only provided GMRES execution details since the computation time of the preconditioner has already been discussed above.
	%\footnote{\color{blue}{Here the comparison between GMRES with SPAI and  GMRES with reusing SPAI  are not given. The reason is, in both cases, the differences are very less, GMRES with  reusing SPAI takes slightly more time (i.e., few seconds) as compare to GMRES with SPAI.}}. 
	%	In the first column, we give the number of AIRGA iterations. The number of times GMRES is executed for each of the linear systems are given in second column. 
	%	The iteration count and computation time of GMRES are given in columns 3 and 4, respectively.
	The structure of this table is same as that of Table \ref{tab:GMRES4K}. 
	As earlier, 
	we notice from this table that solving linear systems by GMRES with SPAI takes less computation time as compared to solving them by GMRES with reusable SPAI. This is again because  of additional matrix-vector products in the reusable SPAI case. Here also, this extra cost is almost negligible when compared to the savings in the preconditioner computation time (as evident in Table \ref{tab:SPAI_SPAI_Update_time_12m}; also see the total GMRES and preconditioner time below).

	%we notice from this table that solving individual linear systems by GMRES with SPAI takes less computation time (i.e., a few minutes) as compared to GMRES with reusing SPAI. This is because, again, when we apply reusing SPAI with GMRES, additional matrix-vector product operation is also performed. 
	%Here, we have not added the computation time of preconditioner. As we have already mentioned earlier, unpreconditioned GMRES does not converge for the given model. %Here, we have taken maximum iteration in GMRES for convergence is $1000$.% and $10^{-06}$ as an stopping tolerance. 
	Table~\ref{tab:gmres_spai_spai_u} 
	%gives the computation time of GMRES with SPAI and GMRES with reusing SPAI. The structure of this table is similar to that of Table \ref{tab:gmres_spai_spai_u4k}.  
	gives the computation time of GMRES plus SPAI (column 2) and GMRES plus reusable SPAI (column 3) at each AIRGA iteration (column 1). As before, it is evident from this table, reusing the SPAI preconditioner leads to about $64 \%$ savings in total time (from $197$ hours $28$ minutes to $72$ hours $06$ minutes).% required for solving all the linear systems.
	%In the first column, we give the number of AIRGA iterations. The admissible computation times are given in columns 2 and 3, respectively. 
	%We can notice from this table that computation time of GMRES with reused SPAI is on an average $60\%$ less than the computation time of GMRES with SPAI.
	% when using basic SPAI preconditioner  and when using SPAI with update preconditioner  for the given model. 
	%We can notice from this table that iterative solves with update take about ? of time as needed for basic iterative solves (or without SPAI update).
	%We can notice from this table that computation time of GMRES with SPAI update is on an average $57\%$ less than the computation time of GMRES with SPAI. 
	%This saving is larger for model size 10000   (we go from around 1 hour to 40 Minutes). Hence, larger the problem more the saving.  
	%We observe that the total computation time is reduced from $197$ hrs $19$ mins to $72$ hrs $01$ mins (approximately $64 \%$) by using  GMRES with reusing SPAI  as compared to GMRES with SPAI. 	

\begin{table*}[!]
	%\tiny
	\centering
	\def\arraystretch{1.0}
			\setlength\tabcolsep{3pt}
	\caption{GMRES computation time for the industrial disk brake model.}
	\label{tab:GMRES}
	\begin{tabular}{|c|c|c|c|c|}
		\hline
		\multirow{2}{*}{\begin{tabular}[c]{@{}c@{}}AIRGA \\ Iterations $(z)$ \end{tabular}}  
		&{{\begin{tabular}[c]{@{}c@{}} No. of \\ Linear \\ Solves \end{tabular}}} & {{\begin{tabular}[c]{@{}c@{}} GMRES \\ Iterations  per \\ Linear Solve \end{tabular}}} & {{\begin{tabular}[c]{@{}c@{}} GMRES Time \\ when Using \\ SPAI \\ (Minutes) \end{tabular}}} & {{\begin{tabular}[c]{@{}c@{}} GMRES Time \\ when Using \\ Reusable SPAI \\ (Minutes) \end{tabular}}} \\ \hline
		1 & 64 & $421$    & $08 $ & $10 $ \\ \hline
		2 & 64 & $426$    & $09 $ & $11 $ \\ \hline
		3 & 64 & $429$    & $10 $ & $12 $  \\ \hline
		4 & 52 & $432$    & $10 $ & $12 $\\ \hline
		\cline{1-5}
		\textbf{Total} & {\bf 244} & \begin{tabular}[c]{@{}c@{}} $64 \times (421 + 426 +  429) $ \\ $+ \ 52 \times 432$
			\\ $\mathbf{= 104,128}$ \end{tabular} &  \begin{tabular}[c]{@{}c@{}} $64 \times (08 +  09 + 10)$ \\ $+\ 52 \times 10$
			\\ $\mathbf{= 2248}$  \end{tabular} &  \begin{tabular}[c]{@{}c@{}} $64 \times (10 +  11 + 12)$ \\ $+\ 52 \times 12$
			\\ $\mathbf{= 2736}$  \end{tabular} \\ \hline
	\end{tabular}
\end{table*}

To demonstrate the quality of the reduced system, we plot the relative $H_2$ error between the transfer function of the original system and the reduced system with respect to the different expansion points (in Figure \ref{fig:acc_analy}). The reduced system considered here is obtained by using GMRES with reusable SPAI. 
	%Figure \ref{fig:acc_analy} gives the relative error between the original and the reduced system for this model.	
These expansion points, denoted by $S$, are computed as $2 \pi f$, where the frequency variable $f$ is linearly spaced between $1$ and $500$. 
	%	are given on the x-axis, and error on the y-axis.
As evident from this figure, the obtained reduced system is good (the error is very small). 
Further, we also observe from this figure that the reduced model is most accurate in 7--10 range of the expansion points. This is because the final expansion points, upon the convergence of the AIRGA algorithm, lie in this range. 
}
\begin{figure}[H]
	\centering
	\includegraphics[width=95mm]{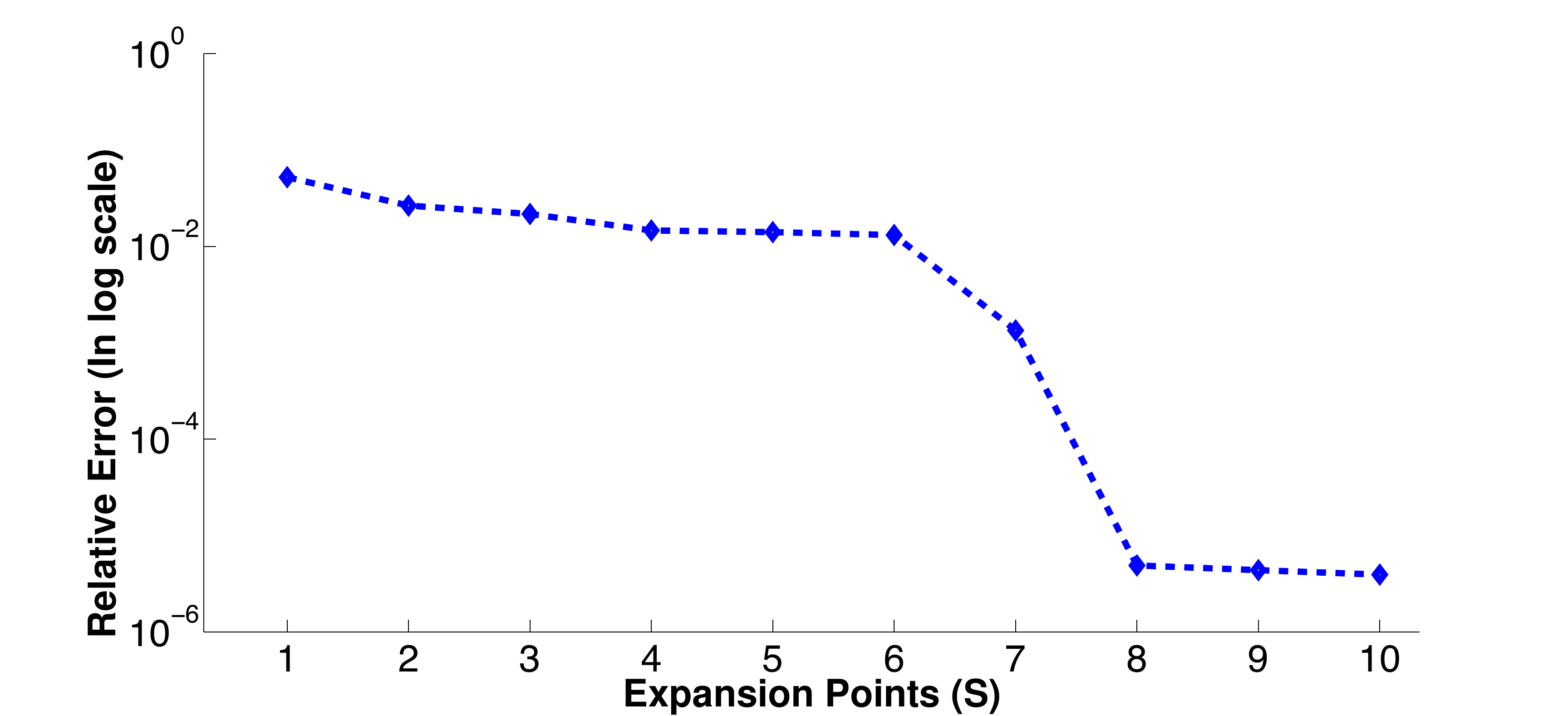}
	\caption{Relative error between the original and reduced system for the industrial disk brake model.}
	\label{fig:acc_analy}       % Give a unique label
\end{figure}

\begin{table}[!]
	%\tiny
	\centering
	\def\arraystretch{1.0}
		\setlength\tabcolsep{3pt}
	\caption{{GMRES with SPAI and reusable SPAI computation time for the industrial disk brake model.}}
	\label{tab:gmres_spai_spai_u}
	\begin{tabular}{|c|c|c|}
		\hline
		{{\begin{tabular}[c]{@{}c@{}}AIRGA \\ Iterations $(z)$ \end{tabular}}} & {{\begin{tabular}[c]{@{}c@{}}GMRES plus\\ SPAI Time  \end{tabular}}} & {{\begin{tabular}[c]{@{}c@{}}GMRES plus\\ Reusable SPAI Time  \end{tabular}}} \\ \hline
		1 & $48 \ \text{hrs} \ 32 \ \text{mins} $ & $23 \ \text{hrs}  \ 40 \ \text{mins}$ \\ \hline
		2 & $49 \ \text{hrs} \ 36 \ \text{mins} $ & $16 \ \text{hrs}  \ 14 \ \text{mins}$    \\ \hline
		3 & $50 \  \text{hrs} \ 40 \ \text{mins}$ & $17 \ \text{hrs}  \ 18 \ \text{mins}$  \\ \hline
		4 & $48 \ \text{hrs} \ 40 \ \text{mins} $& $14 \ \text{hrs}  \ 54 \ \text{mins}$  \\ \hline
		\cline{1-3} \bigstrut
		\ \ \ \ \ \textbf{Total} & $\mathbf{197 \ \text{\bf hrs}  \ 28 \ \text{\bf mins}}$ & $\mathbf{72 \ \text{\bf hrs}  \ 06  \ \text{\bf mins}}$ \\ \hline
	\end{tabular}
\end{table}

%%%%%%%%%%%%%%%%%%%%%%%%%%%%%%%%%%%%%%%%

%\vspace{5cm}
%\newpage
\section[Conclusions \& Future Work]{Conclusions \& Future Work} 
\label{sec:conl_fut}
In  this work, we have focused on MOR of non-parametric dynamical systems, %both linear and bilinear 
specifically on the following three algorithms: AIRGA, BIRKA, and QB-IHOMM. Since solving large and sparse linear systems is a bottleneck in scaling these MOR algorithms for reduction
of large sized dynamical systems, we have proposed reusing of the SPAI preconditioner.

Specifically, we have demonstrated the following: exploitation of the simplicity because of the lack of parameters in reusing preconditioners, multiple ways of reusing preconditioners within the algorithm, efficient implementation to ensure that the savings because of reusing preconditioners are not negated by bad coding, and experimentation on a massively large industrial problem. Numerical experiments show the effectiveness of our approach, where for a problem of size $1.2$ million, we save upto $64 \%$ in the computation time. In absolute terms, this gives  a saving of $5$ days.

In future, we plan to explore two directions. The first is use of the randomized preconditioners in solving linear systems arising in MOR. This is giving promising results. The second is to use the spiking neural networks to optimize the parameters inside the preconditioners.

\appendices
\section{}
\label{app:BIRKA}
{\color{black} In the Algorithm \ref{Algo:BIRKA}, we solve linear systems of equations at lines 4 and 5. We first apply our proposed theory of reusing preconditioners to line 4, which is given as
	\begin{align*}
		vec \left ( V \right ) = \left ( -\Lambda \otimes I_n - I_{r} \otimes K - \sum\limits_{\mathsf{j}=1}^{m} \check{\check{N}}^T_{\mathsf{j}} \otimes  N_{\mathsf{j}}   \right )^{-1}  \\ \left ( \check{\check{F}}^T \otimes F\right ) \ {vec(I_{m})}.
	\end{align*}
	Here, $\Lambda$ is a diagonal matrix comprising of interpolation points, which is updated at the start of the \texttt{while} loop at line 2. %until convergence (e.g., we taken $\Lambda = \Lambda_i$ for $i = 1, \ldots, \mathsf{z}_{max}$).
	Let $A_{z-1} = -\Lambda_{z-1} \otimes I_n - I_{r} \otimes K -  \sum\limits_{\mathsf{j}=1}^{m} \check{\check{N}}^T_{\mathsf{j}} \otimes  N_{\mathsf{j}}$ and $A_{z} = -\Lambda_{z} \otimes I_n - I_{r} \otimes K -  \sum\limits_{\mathsf{j}=1}^{m} \check{\check{N}}^T_{\mathsf{j}} \otimes  N_{\mathsf{j}}$ %(for sake of simplification we assume $\mathsf{j} =1$ here) 
	be the coefficient matrices corresponding to $\Lambda_{z-1}$ and $\Lambda_{z}$, respectively $\left(\text{for} \ z = 1, \ \ldots, \ \mathfrak{z} \ (\text{until covergence})\right)$.  
	% These interpolation points can change at each of the BIRKA iterations. If the difference between $A_{old}$ and $A_{new}$ is small, then one can exploit this while building preconditioners for this sequence of matrices. 
	Expressing $A_{z}$ in terms of $A_{z-1}$, we get
	\begin{align*}
		A_{z} = A_{z-1} \left(I_{nr} + A_{z-1}^{-1} (-\Lambda_{z} \otimes I_n) + A_{z-1}^{-1} (\Lambda_{z-1} \otimes I_n)\right),
	\end{align*}
	where $I_{nr} \in  \mathbb{R}^{n \cdot r \times n \cdot r}$ is the Identity matrix. %Let $P_{old}$ be a good initial preconditioner for $A_{old}$. A reuse of preconditioner can be obtained by enforcing $A_{old} P_{old} \approx A_{new} P_{new}$.
	Now, we enforce 
	\begin{align}\label{eq:hoti_dir}
		A_{z-1}P_{z-1} = A_{z} P_{z} \ \  or
	\end{align}
	\begin{align*}
		A_{z-1}P_{z-1} =  \ A_{z-1} \left(I_{nr} + A_{z-1}^{-1} (-\Lambda_{z} \otimes I_n) + A_{z-1}^{-1} (\Lambda_{z-1} \otimes I_n)\right)  \cdot \\ \left(I_{nr} + A_{z-1}^{-1} (-\Lambda_{z} \otimes I_n) + A_{z-1}^{-1} (\Lambda_{z-1} \otimes I_n)\right)^{-1} P_{z-1} \\
		=  \ A_{z} P_{z}, 
	\end{align*}
	where $P_{z} = \left(I_{nr} + A_{z-1}^{-1} (-\Lambda_{z} \otimes I_n) + A_{z-1}^{-1} (\Lambda_{z-1} \otimes I_n)\right)^{-1} P_{z-1}$.  
	
	Let $Q_{z} = \left(I_{nr} + A_{z-1}^{-1} (-\Lambda_{z} \otimes I_n) + A_{z-1}^{-1} (\Lambda_{z-1} \otimes I_n)\right)^{-1}$, then instead of (\ref{eq:hoti_dir}) we enforce $A_{z-1} P_{z-1} = A_{z} Q_{z} P_{z-1}$. The remaining derivation here is same as earlier (see Section \ref{sec:app-reuse}).   
	%instead of computing a preconditioner by minimizing $\|I - A_{new} P_{new}\|_f^2$ with respect to $P_{new}$, we solve a simpler problem given below    
	% \begin{align*}
	% \min_{Q_{new}} \| A_{old} - A_{new} Q_{new} \|_f^2,
	% \end{align*} 
	%with $P_{new} = Q_{new}P_{old}$. 
	We reuse preconditioners at line 5 similarly.
	% \begin{align*}
	%vec\left ( W \right ) = \left ( -\Lambda \otimes I_n - I_{r} \otimes K^T - \sum\limits_{{\mathsf{j}}=1}^{m}\check{\check{N}}_{\mathsf{j}} \otimes N^T_{\mathsf{j}} \right )^{-1} \left ( \check{\check{C}}^T \otimes C^T\right ) \ {vec(I_{q})}.
	% \end{align*}
	% Again, here also, $\Lambda$ is the interpolation points, which varying inside of the while-loop (line 2) until convergence (e.g., we taken $\Lambda = \Lambda_i$ for $i = 1, \ldots, \mathsf{z}_{max}$). Let $A_1 = -\Lambda_1 \otimes I_n - I_{r} \otimes K^T -  \check{\check{N}}_{\mathsf{j}} \otimes  N^T_{\mathsf{j}} $ and $A_2 = -\Lambda_2 \otimes I_n - I_{r} \otimes K^T -  \check{\check{N}}_{\mathsf{j}} \otimes  N_{\mathsf{j}}^T$ are the coefficient matrices, and rest are same as above. 
	% In Algorithm \ref{Algo:BIRKA-precond}, we propose this reuse of preconditioners theory in the BIRKA algorithm context (Note:  $A_1^{(1)}= -\Lambda^{(1)} \otimes I_n - I_{r} \otimes K -  \check{\check{N}}^T_1 \otimes  N_1$ (which is $A_{old}$), and $A_1^{(z)}= -\Lambda^{(z)} \otimes I_n - I_{r} \otimes K -  \check{\check{N}}^T_1 \otimes  N_1$ (which is $A_{new}$)). Here, parts relevant to solving linear systems are only listed.
} 
 %----------------------------------------------------------------------------
%\newpage
\section{}
\label{app:QB-IHOMM}
{\color{black} In the Algorithm \ref{Algo:QB-IHOMM}, we solve linear systems of equations at line 4 and 10. Again, we first apply our proposed theory of reusing preconditioners to line 4, which is given as
	%	Now, we apply our proposed theory of reuse preconditioners to solve these linear system of equations. The equation given at line 4 is
	\begin{align*}
		X_j(\sigma_i) = [(\sigma_iD-K)^{-1}D]^{j}(\sigma_iD-K)^{-1}F, \ \\  \text{for} \ j = 1, \ \ldots, \ P+Q \ \text{and}  \  i = 1, \ \ldots, \ \ell.
	\end{align*}
	Let $A_{i-1} = \sigma_{i-1}D - K$ and $A_{i} =  \sigma_{i}D - K$ be the two coefficient matrices for different interpolation points $\sigma_{i-1}$ and $\sigma_{i}$, respectively $\left(\text{for} \ i = 1, \ \ldots, \ \ell\right)$. 
	%These interpolation points can change at each QB-IHOMM iterations. Again, if the difference between $A_{old}$  and $A_{new}$ is small, then one can exploit this while building preconditioners for this sequence of matrices.
	Expressing $A_{i}$ in terms of $A_{i-1}$, we get
	\begin{align*}
		A_{i} = A_{i-1}(I  + (\sigma_{i} - \sigma_{i-1})A_{i-1}^{-1}D).
	\end{align*}
	%Let $P_{old}$ be a good initial  preconditioner for $A_{old}$. Here also, a reuse of preconditioner can be obtained  by enforcing $A_{old}P_{old} \approx A_{new}P_{new}$, where $old,new =\{1, \ \ldots, \ \ell \}$ and, as earlier, $\ell$ denotes the number of interpolation points. 
	Now, we enforce 
	\begin{align}\label{eq:hori_dir_QB}
		A_{i-1} P_{i-1} = A_{i} P_{i} \ \  or
	\end{align}
	\begin{align*}
		\begin{split}
			A_{i-1}P_{i-1} =  \ A_{i-1}(I + (\sigma_{i} - \sigma_{i-1})A_{i-1}^{-1}D) 
			\ \cdot \\ (I +  (\sigma_{i}-\sigma_{i-1})A_{i-1}^{-1}D)^{-1}P_{i-1} \\
			= \ A_{i}P_{i},
		\end{split}
	\end{align*}
	where  $P_{i} = (I + (\sigma_{i} - \sigma_{i-1})A_{i-1}^{-1}D)^{-1}P_{i-1}$. 
	
	Let $Q_{i} = (I + (\sigma_{i} - \sigma_{i-1})A_{i-1}^{-1}D)^{-1}$,  then instead of (\ref{eq:hori_dir_QB}) we enforce $A_{i-1}P_{i-1} = A_{i}Q_{i}P_{i-1}$. Again, here also, the remaining derivation is same as earlier (see Section \ref{sec:app-reuse}).
	%instead of computing a preconditioner by minimizing $\|I - A_{new}P_{new}\|_f^2$ with respect to $P_{new}$, we solve a simpler problem given below
	%\begin{equation*}
	%\min_{Q_{new}}\|A_{old}-A_{new}Q_{new}\|^{2}_{f},
	%\end{equation*}
	%with $P_{new} = Q_{new}P_{old}$. 
	We reuse preconditioners at line 10 similarly. 
\section*{Acknowledgment}
We would like to deeply thank Prof. Dr. Heike Fa{\ss}bender (at Institut Computational Mathematics, AG Numerik, Technische Universit{\"a}t Braunschweig, Germany) for discussions and help regarding different aspects of this project.

\begin{IEEEbiography}[{\includegraphics[width=1.05in,height=1.22in,clip]{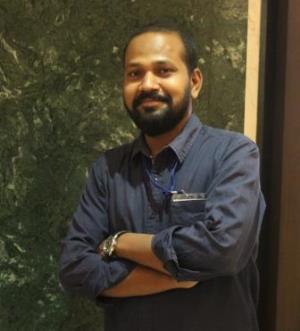}}]{Navneet Pratap Singh}received his bachelor's degree in Computer Science and Engineering from  UPTU, Lucknow, India, and his master's degrees in Modelling and Simulation from the Defence Institute of Advanced Technology, Pune, India. 

He is currently pursuing his Ph.D. degree with IIT Indore. His thesis focuses on Stable Linear Solves with Preconditioner Updates for Model Reduction. His research interests are at the intersection of Computer Science and Mathematics, especially Numerical Linear Algebra, Optimization, Dynamical Systems, and Machine Learning. 

\vspace{-12cm}
\end{IEEEbiography}
\begin{IEEEbiography}[{\includegraphics[width=1.05in,height=1.22in,clip]{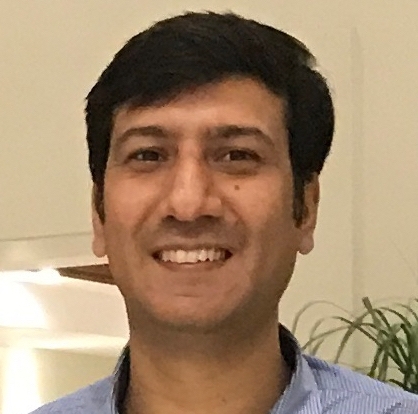}}]{Kapil Ahuja} (B.Tech.: IIT (BHU), India; M.S. and Ph.D.: Virginia Tech, USA; Postdoctoral Research Fellow: Max Planck Institute, Germany) has a varied background, including degrees in Computer Science, Mathematics, and Mechanical Engineering. 
	
He is currently an Associate Professor in Computer Science and Engineering at IIT Indore (India). In the past, he has been a visiting professor at TU Braunschweig (Germany), TU Dresden (Germany), and Sandia National Labs (USA). Dr. Ahuja is working on Mathematics of Data Science as well as Computational Science. Specifically, Artificial Intelligence, Machine Learning, Numerical Methods, and Optimisation.
\end{IEEEbiography}
\EOD
\end{document}